\documentclass[10pt,draft,reqno]{amsart}
     \makeatletter
     \def\section{\@startsection{section}{1}%
     \z@{.7\linespacing\@plus\linespacing}{.5\linespacing}%
     {\bfseries
     \centering
     }}
     \def\@secnumfont{\bfseries}
     \makeatother
\newtheorem{theorem}{Theorem}[section]
\newtheorem{lemma}[theorem]{Lemma}
\newtheorem{proposition}[theorem]{Proposition}

\theoremstyle{definition}
\newtheorem{definition}[theorem]{Definition}

\theoremstyle{remark}
\newtheorem{remark}[theorem]{Remark}
\numberwithin{equation}{section}
\def\real{{\mathord{{\rm I\kern-2.8pt R}}}}
\def\inte{{\mathord{{\rm I\kern-2.8pt N}}}}

\def\sZZ{{\rm Z\kern-2.8ptem{}Z}}

\def\z{{\mathchoice
{\sZZ}
{\sZZ}
{\rm Z\kern-0.30em{}Z}
{\rm Z\kern-0.25em{}Z} }}
\def\sQQ{{\kern 0.27em \vrule height1.45ex width0.03em depth0em
\kern-0.30em \rm Q}}
\def\qu{{\mathchoice
{\sQQ}
{\sQQ}
{\kern 0.225em \vrule height1.05ex width0.025em depth0em \kern-0.25em \rm Q}
{\kern 0.180em \vrule height0.78ex width0.020em depth0em \kern-0.20em \rm Q}
}}
\def\sCC{{\kern 0.27em \vrule height1.45ex width0.03em depth0em
\kern-0.30em \rm C}}
\def\complex{{\mathchoice
{\sCC}
{\sCC}
{\kern 0.225em \vrule height1.05ex width0.025em depth0em \kern-0.25em \rm C}
{\kern 0.180em \vrule height0.78ex width0.020em depth0em \kern-0.20em \rm C}
}}

\newcommand{\ba}{\begin{array}}
\newcommand{\ea}{\end{array}}
\newcommand{\be}{\begin{equation}}
\newcommand{\ee}{\end{equation}}
\newcommand{\bea}{\begin{eqnarray}}
\newcommand{\eea}{\end{eqnarray}}
\newcommand{\beaa}{\begin{eqnarray*}}
\newcommand{\eeaa}{\end{eqnarray*}}

\def\z{\zeta}

\font\tenmath=msbm10 \font\sevenmath=msbm7 \font\fivemath=msbm5
\newfam\mathfam \textfont\mathfam=\tenmath
\scriptfont\mathfam=\sevenmath \scriptscriptfont\mathfam=\fivemath

\def \={{\buildrel {\rm (law)} \over =}}

\def\qed{ \hfill \vrule width.25cm height.25cm depth0cm\smallskip}

\newcommand{\basa}{\begin{assumption}}
\newcommand{\easa}{\end{assumption}}

\newcommand{\bas}{\begin{assum}}
\newcommand{\eas}{\end{assum}}

\begin{document}
\title[short title for running heading]{Self-similarity parameter estimation and reproduction property for
non-Gaussian Hermite processes}
\author[Alexandra Chronopoulou]{Alexandra Chronopoulou}
\thanks{* This author's research partially supported by NSF grant DMS 0606615}
\address{Alexandra Chronopoulou: Dept. Statistics, Purdue University, West Lafayette, Indiana IN 47906, USA}
\email{achronop@stat.purdue.edu}

\author[Frederi G. Viens]{Frederi G. Viens}
\thanks{* This author's research partially supported by NSF grants DMS 0606615 and 0907321}
\address{Frederi G. Viens: Department of Statistics, Purdue University, West Lafayette, Indiana IN 47906, USA}
\email{viens@stat.purdue.edu}
\author[Ciprian A. Tudor]{Ciprian A. Tudor}
\thanks{* Associate member of the team SAMM, Universit\'e de Panth\'eon-Sorbonne Paris 1 }
\address{Ciprian A. Tudor: Laboratoire Paul Painlev\'e, Universit\'e de Lille 1, F-59655 Villeneuve d'Ascq, France}
\email{tudor@math.univ-lille1.fr}


\begin{abstract}Let $(Z^{(q,H)}_{t})_{t\in [0,1]}$ be a Hermite processes of order $q$ and  with Hurst parameter $H\in
( \frac{1}{2}, 1)$. This process is $H$-self-similar, it has stationary increments and it exhibits long-range dependence. This
class contains the fractional Brownian motion (for $q=1$) and the Rosenblatt process (for $q=2$). We study in this paper the
variations of $Z^{(q,H)}$ by using multiple Wiener -It\^o stochastic integrals and Malliavin calculus. We prove
a reproduction property for this class of processes in the sense that the terms appearing in the chaotic decomposition of the
their variations give birth to other Hermite processes of different orders and with different Hurst parameters. We apply our
results to construct a consistent estimator for the self-similarity parameter from discrete observations of a Hermite process.
\end{abstract}
\maketitle
\section{Introduction}
\subsection{Background and motivation}
The variations of a stochastic process play a crucial role in its
probabilistic and statistical analysis. Best known is the quadratic
variation of a semimartingale, whose role is crucial in It\^{o}'s formula
for semimartingales; this example also has a direct utility in practice, in
estimating unknown parameters, such as volatility in financial models, in
the so-called \textquotedblleft historical\textquotedblright\ context. For
self-similar stochastic processes the study of their variations constitutes
a fundamental tool to construct good estimators for the self-similarity
parameter. These processes are well suited to model various phenomena where long
memory is an important factor (internet traffic, hydrology, econometrics,
among others). The most important modeling task is then to determine or
estimate the self-similarity parameter, because it is also typically responsible
for the process's long memory and other regularity properties. Consequently,
estimating this parameter represents an important research direction in theory
and practice. Several approaches, such as wavelets, variations, maximum
likelihood methods, have been proposed. We refer to the monograph \cite%
{Beran} for a complete exposition.
The family of Gaussian processes known as fractional Brownian motion (fBm)
is particularly interesting, and most popular among self-similar processes,
because of fBm's stationary increments, its clear similarities and
differences with standard Brownian motion, and the fact that its
self-similarity parameter $H$, known as the Hurst parameter, is also clearly
interpreted as the memory length parameter (the correlation of unit length
increments $n$ time units apart decays slowly at the speed $n^{2H-2}$), and
the regularity parameter (fBm is $\alpha$-H\"{o}lder continuous on any
bounded time interval for any $\alpha<H$).
Soon after fBm's inception, the study of its variations began in the 1970's
and 1980's; of interest to us in the present article are several such
studies of variations which uncovered a generalization of fBm to
non-Gaussian processes known as the Rosenblatt process and other Hermite
processes: \cite{BrMa}, \cite{DM}, \cite{GS}, \cite{Ta1} or \cite{Ta2}. We
briefly recall some relevant basic facts. We consider $(B_{t}^{H})_{t\in
\lbrack 0,1]}$ a fractional Brownian motion with Hurst parameter $H\in (0,1)$%
. As such, $B^{H}$ is the continuous centered Gaussian process with  covariance function
\begin{equation*}
R_{H}(t,s)=\mathbf{E}\left[ B^{H}_{t}B^{H}_{s}\right] =\frac{1}{2}\big(%
s^{2H}+t^{2H}-|t-s|^{2H}\big),\quad s,t\in \lbrack 0,1].
\end{equation*}%
Equivalently, $B^{H}_{0}=0$ and $\mathbf{E}\left[ \left( B^{H}_{s}-B^{H}_{t}\right) ^{2}%
\right] =\left\vert t-s\right\vert ^{2H}$. It is the only Gaussian process $H
$ which is self-similar with stationary increments. Consider $%
0=t_{0}<t_{1}<\ldots <t_{N}=1$ a partition of the interval $[0,1]$ with $%
t_{i}=\frac{i}{N}$ for $i=0,\ldots, N$ and define the following $q$ variations
\begin{equation*}
v_{N}^{(q)}=\sum_{i=0}^{N-1}H_{q}\left( \frac{B_{t_{i+1}}^{H}-B_{t_{i}}^{H}}{%
\left( \mathbf{E}\left[ \left( B_{t_{i+1}}^{H}-B_{t_{i}}^{H}\right) ^{2}%
\right] \right) ^{\frac{1}{2}}}\right) =\sum_{i=0}^{N-1}H_{q}\left( \frac{\left( B_{t_{i+1}}^{H}-B_{t_{i}}^{H}\right)}{\left(
t_{i+1}-t_{i}\right) ^{H}}
\right)
\end{equation*}%
where $H_{q}$ with $q\geq 2$ represents the Hermite polynomial of degree $q$. Then it follows from the above references that for:
\begin{itemize}
\item for $0\leq H < 1-\frac{1}{2q}$  the limit in distribution of $N^{-%
\frac{1}{2}}v_{N}^{(q)}$ is a centered Gaussian random variable,
\item for $H=1-\frac{1}{2q}$  the limit in distribution of $(N \log N)^{-%
\frac{1}{2}}v_{N}^{(q)}$ is a centered Gaussian random variable,
\item for $1-\frac{1}{2q}<H<1$ the limit of $N^{q(1-H)-1}v_{N}^{(q)}$ is a
Hermite random variable of order $q$ with self-similarity parameter $q(H-1)+1.$
\end{itemize}
This latter random variable is non-Gaussian; it equals the value at time $1$ of a Hermite process, which is a stochastic
process in the $q$th Wiener chaos with the same covariance structure as fBm; as such, it has stationary increments and shares
the same self-similarity, regularity, and long memory properties as fBm; see Definition \ref{defHermite}. We also mention that
very recently, various interesting results have been proven for weighted power variations of stochastic processes such as
fractional Brownian motion (see \cite{No}), fractional Brownian sheets (see \cite{Rev}),  iterated Brownian motion (see
\cite{NoPe}) or the solution of the stochastic heat equation (see \cite{Sw} or \cite{BuSw}).
Because of a natural coupling, the last limit above also holds in $%
L^{2}(\Omega )$ (see \cite{NNT}). In the critical case $H=1-\frac{1}{2q}$
the limit is still Gaussian but the normalization involves a logarithm.
These results are widely applied to estimation problems; to avoid the
barrier $H=\frac{3}{4}$ that occurs in the case $q=2$, one can use
\textquotedblleft higher-order filters\textquotedblright , which means that
the increments of the fBm are replaced by higher-order increments, such as
discrete versions of higher-order derivatives, in order to obtain a Gaussian
limit for any $H$ (see \cite{LaIs}, \cite{GuLe}, \cite{coeur}).
The appearance of Hermite random variables in the above non-central limit
theorems begs the study of Hermite processes as such. Their practical aspects
are striking: they provide a wide class of processes from which to model
long memory, self-similarity, and H\"{o}lder-regularity, allowing
significant deviation from fBm and other Gaussian processes, without having
to invoke non-linear stochastic differential equations based on fBm, and the
notorious issues associated with them (see \cite{N}). Just as in the case of
fBm, the estimation of the Hermite process's parameter $H$ is crucial for proper
modeling; to our knowledge it has not been treated in the literature. We
choose to tackle this issue by using variations methods, to find out how the
above central and non-central limit theorems fit in a larger picture.
\subsection{Main results: summary and discussion}
In this article, we show results that are interesting from a theoretical
viewpoint, such as the reproduction properties (the variations of Hermite
processes give birth to other Hermite processes); and we provide an
application to parameter estimation, in which care is taken to show that the
estimators can be evaluated practically.
Let $Z^{(q,H)}$ be a Hermite process of order $q$ with selfsimilarity parameter $H\in (\frac{1}{2}, 1)$ as defined in Definition
\ref{defHermite}.
Define the \emph{centered quadratic variation statistic}%
\begin{equation}\label{VN}
V_{N}:=\frac{1}{N}\sum_{i=0}^{N-1}\left[ \frac{\left( Z_{\frac{i+1}{N}%
}^{(q,H)}-Z_{\frac{i}{N}}^{(q,H)}\right) ^{2}}{\mathbf{E}\left[ \left( Z_{%
\frac{i+1}{N}}^{(q,H)}-Z_{\frac{i}{N}}^{(q,H)}\right) ^{2}\right] }-1\right].
\end{equation}%
Also for $H\in (1/2,1)$, and $q\in \mathbf{N}\setminus \{0\}$, we define a constant
which will recur throughout this article:%
\begin{eqnarray}
d(H, q) &:&=\frac{\left( H(2H-1)\right) ^{1/2}}{\left( q!(H'(2H'-1)) ^{q}\right) ^{1/2}}, \hskip0.2cm  H^{\prime }
=1+\frac{\left( H-1\right) }{q}.  \label{dH}
\end{eqnarray}%
We prove that, under suitable normalization, this sequence converges in $%
L^{2}(\Omega )$ to a Rosenblatt random variable.
\begin{theorem}
\label{thmgen1}Let $H\in(1/2,1)$ and $q\in\mathbf{N}\setminus \{0\}$. Let $Z^{(q,H)} $ be a Hermite process of order $q$ and
self-similarity parameter $H$ (see
Definition \ref{defHermite}). Let $c$ be an explicit positive constant (it is defined in Proposition \ref{PropLim}).
Then $cN^{(2-2H)/q}V_{N}$ converges in $L^{2}\left( \Omega\right) $ to a
standard Rosenblatt random variable with parameter $H^{\prime\prime}:=\frac{%
2(H-1)}{q}+1$, that is, to the value at time $1$ of a Hermite process of
order 2 and self-similarity parameter $H^{\prime\prime}$.
\end{theorem}
The Rosenblatt random variable is a double integral with respect to the same Wiener process used to define the Hermite
process; it is thus an element of the second Wiener chaos. Our result shows that fBm is the only Hermite process for which
there exists a range of parameters allowing normal convergence of the quadratic variation, while for all other Hermite
processes, convergence to a second chaos random variable is universal. Our proofs are based on chaos expansions into multiple
Wiener integrals and Malliavin calculus. The main line of the proof is as follows: since the variable $Z_{t}^{(q,H)}$ is an
element of the $q$th Wiener chaos, the product formula for multiple integrals implies that the statistics $V_{N}$ can be
decomposed into a sum of multiple integrals from the order $2$ to the order $2q$. The dominant term in this decomposition,
which gives the final renormalization order $N^{(2-2H)/q}$, is the term which is a double Wiener integral, and one proves it
\emph{always }converges to a Rosenblatt random variable; all other terms are of much lower orders, which is why the only
remaining term, after renormalization, converges to a second chaos random variable. The difference with the fBm case comes
from the limit of the term of order 2, which in that case is sometimes Gaussian and sometimes Rosenblatt-distributed,
depending on the value of $H$.
We also study the limits of the other terms in the decomposition of $V_{N}$,
those of order higher than $2$, and we obtain interesting facts: all these
terms, except the term of highest order $2q$, have limits which are Hermite
random variables of various orders and self-similarity parameters. We call this
\emph{the reproduction property} for Hermite processes, because from one
Hermite process of order $q$, one can reconstruct other Hermite processes of
all lower orders. The exception to this rule is that the normalized term of
highest order $2q$ converges to a Hermite r.v. of order $2q$ if $H>3/4$, but
converges to Gaussian limit if $H\in(1/2,3/4]$. Summarizing, we have the
following.
\begin{theorem}
\label{thmgen2}Let $Z^{(q,H)}$ be again a Hermite process, as in the
previous theorem. Let $T_{2n}$ be the term of order $2n$ in the Wiener chaos
expansion of $V_{N}$: this is a multiple Wiener integral of order $2n$, and
we write $V_{N}=\sum_{n=1}^{q}c_{2n}T_{2n}\ $where $c_{2q-2k}=k!\binom{q}{k}%
^{2}$.
\begin{itemize}
\item For every $H\in(1/2,1)$ and for every $k=1,\ldots,q-1$ we have
convergence of the expression $\left( z_{k,H}\right) ^{-1}N^{(2-2H^{\prime })k}T_{2k}$ in $L^{2}\left( \Omega\right) $ to
$Z^{(r,K)}$, a Hermite random variable of order $r=2k$ with self-similarity parameter $K=2k(H^{\prime}-1)+1 $, where $z_{k,H}$
is a constant.
\item For every $H\in(1/2,3/4)$, with $k=q$, we have convergence of $%
x_{1,H, q}^{-1/2}\sqrt{N}T_{2q}$ to a standard normal distribution, with $%
x_{1,H, q}$  a positive  constant depending on $H$ and $q$.
\item For every $H\in(3/4,1)$, with $k=q$, we have convergence of $%
x_{2,H, q}^{-1/2}$ $N^{2-2H}$ $T_{2q}$ in $L^{2}(\Omega)$ to $Z^{2q,2H-1}$, a Hermite
r.v. of order $2q$ with parameter $2H-1$; with $ x_{2,H,q}$ a positive constant depending on $H$ and $q$.
\item For $H=3/4$, with $k=q$, we have convergence of $\sqrt{N/\log N}%
x_{3,H, q}^{-1/2}T_{2q}$ to a standard normal distribution, with $x_{3,H}$ a positive constant depending on $H$ and $q$.
\end{itemize}
\end{theorem}

Some of the aspects of this theorem had been discovered in the case of $q=2$
(Rosenblatt process) in \cite{TV}. In that paper, the existence of a
higher-chaos-order term with normal convergence had been discovered for the
Rosenblatt process with $H<3/4$, while the case of $H\geq3/4$ had not been
studied. The entire spectrum of convergences in Theorem \ref{thmgen2} was
not apparent in \cite{TV}, however, because it was unclear whether the term $%
T_{2}$'s convergence to a Rosenblatt r.v. was due to the fact that we were
dealing with input coming from a Rosenblatt process, or whether it was a
more general function of the structure of the second Wiener chaos; here we
see that the second alternative is true.
The paper \cite{TV} also exhibited a remarkable structure of the Rosenblatt data when $%
H<3/4$. In that case, as we see in Theorem \ref{thmgen2}, there are only two
terms in the expansion of $V_{N}$, $T_{2}$ and $T_{4}$; moreover, and this
is the remarkable feature, the proper normalization of the term $T_{2}$
converges to none other than the Rosenblatt r.v. at time $1$. Since this
value is part of the observed data, one can subtract it to take advantage of
the Gaussian limit of the renormalized $T_{4}$, including an application to
parameter estimation in \cite{TV}. In Theorem \ref{thmgen2} above, if $q>2$,
even if $H\leq3/4$, by which a Gaussian limit can be constructed from the
renormalized $T_{2q}$, we still have at least two other terms $%
T_{2},T_{4},\cdots T_{2q-2}$, and all but at most one of these will converge in $L^{2}\left( \Omega\right) $ to Hermite
processes with different orders, all different from $q$, which implies that they are not directly observed. This means our
Theorem \ref{thmgen2} proves that the operation performed with Rosenblatt data, subtracting an observed quantity to isolate
$T_{2q}$ and its Gaussian asymptotics, is not possible with any higher-order Hermite processes. The last aspect of this paper
applies Theorem \ref{thmgen1} to estimating the parameter $H$. Let $S_{N}$ be the empirical mean of the individual squared
increments
\begin{equation*}
S_{N}=\frac{1}{N}\sum_{i=0}^{N-1}\left( Z_{\frac{i+1}{N}}^{(q,H)}-Z_{\frac {i%
}{N}}^{(q,H)}\right) ^{2}
\end{equation*}
\ \ and let
\begin{equation*}
\hat{H}_{N}=\left( \log S_{N}\right) /\left( 2\log N\right) .
\end{equation*}
We show that $\hat{H}_{N}$ is a strongly consistent estimator of $H$, and we
show asymptotic Rosenblatt distribution for $N^{2\left( 1-H\right) /q}\left(
H-\hat{H}\right) \log N$. The fact that this estimator fails to be
asymptotically normal is not a problem in itself. What is more problematic
is the fact that if one tries to check ones assumptions on the data by
comparing the asymptotics of $\hat{H}$ with a Rosenblatt distribution, one
has to know something about the normalization constant $N^{2\left(
1-H\right) /q}$. Here, we prove in addition that one may replace $H$ in this
formula by $\hat{H}_{N}$, so that the asymptotic properties of $\hat{H}_{N}$
can actually be checked.
\begin{theorem}
\label{thmgen3}The estimator $\hat{H}_{N}$ is strongly consistent, i.e. $%
\lim_{N\rightarrow\infty}\hat{H}_{N}=H$ almost surely. Moreover there exists
a standard Rosenblatt random variable $R$ with self-similarity parameter $%
1+2(H-1)/q$ such that%
\begin{equation*}
\lim_{N\rightarrow\infty}\mathbf{E}\left[ \left\vert 2N^{2\left( 1-\hat {H}%
_{N}\right) /q}\left( H-\hat{H}_{N}\right) \log
N-c_{2}c_{1,H}^{1/2}R\right\vert \right] =0,
\end{equation*}
Here $a\left( H^{\prime}\right)=\left( 1+(H-1)/q\right) \left( 1+2(H-1)/q\right) .$
\end{theorem}
Replacing the constant $c_{1,H}$ by its value in terms of $\hat{H}_{N}$ instead of $H$ also seems to lead to the above
convergence. However, this article does not present any numerical results illustrating model validation based on the above
theorem; moreover such applications would be much more sensitive to the convergence speed than to the actual constants;
therefore we omit the proof of this further improvement on $c_{1,H}$.

\section{Preliminaries\label{PrelimSec}}
\subsection{Multiplication in Wiener Chaos}
Let $(W_{t})_{t\in \lbrack 0,1]}$ be a classical Wiener process on a
standard Wiener space $\left( \Omega ,{\mathcal{F}},\mathbf{P}\right) $. If $%
f\in L^{2}([0,1]^{m})$ with $m\geq 1$ integer, we introduce the multiple
Wiener-It\^{o} integral of $f$ with respect to $W$. We refer to \cite{N} for
a detailed exposition of the construction and the properties of multiple
Wiener-It\^{o} integrals.
Let $f\in{\mathcal{S}}$, i.e. $f$ is an elementary function, meaning that%
\begin{equation*}
f=\sum_{i_{1},\ldots,i_{m}}c_{i_{1},\ldots i_{m}}1_{A_{i_{i}}\times
\ldots\times A_{i_{m}}}
\end{equation*}
where $i_{1},\ldots,i_{m}$ describes a finite set and the coefficients satisfy $c_{i_{1},\ldots i_{m}}=0$ if two indices $%
i_{k}$ and $i_{l}$ are equal and the sets $A_{i}\in{\mathcal{B}}([0,1])$ are
disjoints.
For such a step function $f$ we define
\begin{equation*}
I_{m}(f)=\sum_{i_{1}, \ldots, i_{m}}c_{i_{1}, \ldots i_{m}
}W(A_{i_{1}})\ldots W(A_{i_{m}})
\end{equation*}
where we put $W([a,b])= W_{b}-W_{a}$.
It can be seen that the application $I_{m}$ constructed above from ${%
\mathcal{S}}$ to $L^{2}(\Omega)$ satisfies
\begin{equation}
\mathbf{E}\left[ I_{n}(f)I_{m}(g)\right] =n!\langle f,g\rangle
_{L^{2}([0,1]^{n})}\mbox{ if }m=n  \label{isom}
\end{equation}%
and
\begin{equation*}
\mathbf{E}\left[ I_{n}(f)I_{m}(g)\right] =0\mbox{ if }m\not=n.
\end{equation*}%
It also holds that $I_{n}(f)=I_{n}\left( \tilde{f}\right)$
where $\tilde{f}$ denotes the symmetrization of $f$ defined by $\tilde{f}%
(x_{1},\ldots ,x_{n})=\frac{1}{n!}\sum_{\sigma \in _{n}}f(x_{\sigma (1)},\ldots ,x_{\sigma (n)}).$
Since the set ${\mathcal{S}}$ is dense in $L^{2}([0,1]^{n})$ the mapping $%
I_{n}$ can be extended to a linear continuous operator from from $L^{2}([0,1]^{n})$ to $%
L^{2}(\Omega)$ and the above properties hold true for this extension. Note
also that $I_{n}(f)$ with $f$ symmetric  can be viewed as an iterated stochastic integral
\begin{equation*}
I_{n}(f)=n!\int_{0}^{1}\int_{0}^{t_{n}}\ldots\int_{0}^{t_{2}}f(t_{1},%
\ldots,t_{n})dW_{t_{1}}\ldots dW_{t_{n}};
\end{equation*}
here the integrals are of It\^{o} type; this formula is easy to show for
elementary $f$'s, and follows for general symmetric function $f\in L^{2}([0,1]^{n})$ by a
density argument.
We recall the product for two multiple integrals (see \cite{N}): if $f\in
L^{2}([0,1]^{n})$ and $g\in L^{2}([0,1]^{m})$ are symmetric then it holds that
\begin{equation}
I_{n}(f)I_{m}(g)=\sum_{l=0}^{m\wedge n}l!\binom{m}{l}\binom{n}{l}%
I_{m+n-2l}(f\otimes _{l}g)  \label{prod}
\end{equation}%
where the contraction $f\otimes _{l}g$ belongs to $L^{2}([0,1]^{m+n-2l})$
for $l=0,1,\ldots ,m\wedge n$ and it is given by
\begin{eqnarray}
&&(f\otimes _{l}g)(s_{1},\ldots ,s_{n-l},t_{1},\ldots ,t_{m-l}) = \notag \\
&&\int_{[0,1]^{l}}f(s_{1},\ldots ,s_{n-l},u_{1},\ldots ,u_{l})g(t_{1},\ldots ,t_{m-l},u_{1},\ldots ,u_{l})du_{1}\ldots du_{l}.
\label{contra}
\end{eqnarray}
\subsection{The Hermite process}
Recall that the fractional Brownian motion process $(B_{t}^{H})_{t\in \lbrack 0,1]}$
with Hurst parameter $H\in (0,1)$ can be written as
\begin{equation*}
B_{t}^{H}=\int_{0}^{t}K^{H}(t,s)dW_{s},\quad t\in \lbrack 0,1]
\end{equation*}%
where $(W_{t},t\in \lbrack 0,T])$ is a standard Wiener process,  $%
K^{H}\left( t,s\right)  =c_{H}s^{\frac{1}{2}%
-H}\int_{s}^{t}(u-s)^{H-\frac{3}{2}}u^{H-\frac{1}{2}}du$ if $t>s$ (and  it is zero otherwise), $%
c_{H}=\left( \frac{H(2H-1)}{\beta (2-2H,H-\frac{1}{2})}\right) ^{\frac{1}{2}}
$ and $\beta (\cdot ,\cdot )$ is the Beta function. For $t>s$, the kernel's
derivative is $\frac{\partial K^{H}}{\partial t}(t,s)=c_{H}\left( \frac{s}{t}%
\right) ^{\frac{1}{2}-H}(t-s)^{H-\frac{3}{2}}$. Fortunately we will not need
to use these expressions explicitly, since they will be involved below only
in integrals whose expressions are known.

We will denote by $(Z_{t}^{(q,H)})_{t\in\lbrack0,1]}$ the Hermite process \textit{with self-similarity parameter }$H\in\left(
1/2,1\right) $. Here $q\geq1$ is an integer. The Hermite process can be defined in two ways: as a multiple
integral with respect to the standard Wiener process $(W_{t})_{t\in%
\lbrack0,1]}$; or as a multiple integral with respect to a fractional Brownian motion with suitable Hurst parameter. We adopt the first approach throughout the paper. We refer  to \cite{NNT} of \cite{T} for the
following integral representation of Hermite processes.

\begin{definition}
\label{defHermite} The Hermite process $(Z^{(q,H)}_{t})_{t\in[0,1]}$ of order $%
q\geq1$ and with self-similarity parameter $H\in(\frac{1}{2}, 1)$ is given by
\begin{equation}
Z_{t}^{(q,H)}=d(H,q)\int_{[0,t]^q}\hspace*{-0.05in}dW_{y_{1}}\ldots
dW_{y_{q}}\left( \int_{y_{1}\vee\ldots\vee
y_{q}}^{t}\partial_{1}K^{H^{\prime}}(u,y_{1})\ldots\partial_{1}K^{H^{%
\prime}}(u,y_{q})du\right)\label{z1}
\end{equation}
for $t\in \lbrack0,1]$,  where $K^{H^{\prime}}$ is the usual kernel of the fractional Brownian motion
and
\begin{equation}
H^{\prime}=1+\frac{H-1}{q}\Longleftrightarrow(2H^{\prime}-2)q=2H-2.
\label{H'}
\end{equation}
\end{definition}
Of fundamental importance is the fact that the covariance of $Z^{\left(
q,H\right) }$ is identical to that of fBm, namely
\begin{equation*}
\mathbf{E}\left[ Z^{\left( q,H\right) }_{s}Z^{\left( q,H\right) }_{t} \right] =\frac{1}{2}(t^{2H}+s^{2H}-|t-s|^{2H}).
\end{equation*}%
The constant $d(H,q)$, given in (\ref{dH}) on page \pageref{dH}, is chosen to avoid any additional multiplicative constants.
We stress that $Z^{\left( q,H\right) }$ is far from Gaussian for $q>1$, since it is formed of multiple Wiener integrals of
order $q$. \vskip0.3cm
 The basic properties of the Hermite process are listed below: i) the Hermite process $Z^{(q,H)}$ is
$H$-selfsimilar and it has stationary increments; ii) the mean square of the increment is given by
\begin{equation}
\mathbf{E}\left[ \left\vert Z_{t}^{(q,H)}-Z_{s}^{(q,H)}\right\vert ^{2}%
\right] =|t-s|^{2H};  \label{canonmetric}
\end{equation}%
as a consequence, it follows from Kolmogorov's continuity criterion that $Z^{(q,H)}$
has H\"{o}lder-continuous paths of any order $\delta <H$; iii)  it exhibits long-range dependence in the sense that $\sum_{n\geq 1}%
\mathbf{E}\left[ Z_{1}^{(q,H)}(Z_{n+1}^{(q,H)}-Z_{n}^{(q,H)})\right] =\infty.$
 In fact, the summand in this series is of order $n^{2H-2}$. This property
is identical to that of fBm since the processes share the same covariance structure, and the property is well-known for fBm
with $H>1/2$; iv) if $q=1$ then $Z^{(1,H)}$ is a standard Brownian motion with Hurst parameter $H$ while for $q\geq2$ this
stochastic process is not Gaussian. In the case $q=2$ this stochastic process is known as \emph{the Rosenblatt process}.

\section{Variations of the Hermite process\label{VarSec}}

Since $\mathbf{E}\left( Z_{\frac{i+1}{N}}^{(q,H)}-Z_{\frac{i}{N}}^{(q,H)}\right) ^{2}=N^{-2H}$ and by (\ref{canonmetric}), the
centered quadratic variation statistic $V_{N}$ given in the introduction can be written as
\begin{equation*}
V_{N}=N^{2H-1}\sum_{i=0}^{N-1}\left[ \left( Z_{\frac{i+1}{N}}^{(q,H)}-Z_{\frac{i%
}{N}}^{(q,H)}\right) ^{2}-N^{-2H}\right] .
\end{equation*}
Let $I_{i}:=[\frac{i}{N},\frac{i+1}{N}]$. In preparation for calculating the variance of $V_{N}$ we will find an explicit
expansion of $V_{N}$ in Wiener chaos. We have $Z_{\frac{i+1}{N}}^{(q,H)}-Z_{\frac{i}{N}}^{(q,H)}=I_{q}\left( f_{i,N}\right)$
where we denoted by $f_{i,N}(y_{1},\ldots ,y_{q})$
the expression
\begin{eqnarray}
&&1_{[0,\frac{i+1}{N}]}(y_{1}\vee \ldots \vee y_{q})d(H,q)\int_{y_{1}\vee
\ldots \vee y_{q}}^{\frac{i+1}{N}}\partial _{1}K^{H^{\prime
}}(u,y_{1})\ldots \partial _{1}K^{H^{\prime }}(u,y_{q})du  \notag \\
&&-1_{[0,\frac{i}{N}]}(y_{1}\vee \ldots \vee y_{q})d(H,q)\int_{y_{1}\vee
\ldots \vee y_{q}}^{\frac{i}{N}}\partial _{1}K^{H^{\prime }}(u,y_{1})\ldots
\partial _{1}K^{H^{\prime }}(u,y_{q})du.  \notag
\end{eqnarray}%
Using the product formula for multiple integrals (\ref{prod}), we obtain
\begin{equation}\notag
I_{q}(f_{i,N})I_{q}(f_{i,N})=\sum_{l=0}^{q}l!\binom{q}{k}^{2}I_{2q-2l}\left( f_{i,N}\otimes _{l}f_{i,N}\right)
\end{equation}
where the $f\otimes _{l}g$ denotes the $l$-contraction of the functions $f$
and $g$ given by (\ref{contra}). Let us compute these contractions; for $l=q$
we have
\begin{equation*}
\left\langle f_{i,N}\otimes _{q}f_{i,N} \right\rangle=q!\langle f_{i,N},f_{i,N}\rangle
_{L^{2}([0,1])^{\otimes q}}=\mathbf{E}\left[ \left( Z_{\frac{i+1}{N}%
}^{(q,H)}-Z_{\frac{i}{N}}^{(q,H)}\right) ^{2}\right] =N^{-2H}.
\end{equation*}
Throughout the paper the notation $\partial_{1}K(t,s) $ will be used for $ \partial_{1}K^{H'}(t,s)$. For $l=0$ we have
\begin{eqnarray*}
\left\langle f_{i,N}\otimes _{0}f_{i,N})(y_{1},\ldots y_{q},z_{1},\ldots ,z_{q}\right\rangle
&=&(f_{i,N}\otimes f_{i,N})(y_{1},\ldots y_{q},z_{1},\ldots ,z_{q}) \\
&=&f_{i,N}(y_{1},\ldots ,y_{q})f_{i,N}(z_{1},\ldots ,z_{q})
\end{eqnarray*}
while for $1\leq k\leq q-1$
\begin{eqnarray*}
&&\left\langle f_{i,N}\otimes _{k}f_{i,N})(y_{1},\ldots y_{q-k},z_{1},\ldots
,z_{q-k}\right\rangle =d(H,q)^{2}\int_{[0,1]^{k}}d\alpha _{1}\ldots d\alpha
_{k} \\
&&\hspace*{-0.2in} \hspace*{-0.15in}\left( \mathbf{1}_{i+1,q-k}^{y_{i}}\mathbf{1}_{i+1,k}^{\alpha
_{i}}\int_{\mathbf{I}_{i+1,k}^{y}}du\partial _{1}K(u,y_{1})\ldots \partial
_{1}K(u,y_{q-k})\partial _{1}K(u,\alpha _{1})\ldots \partial _{1}K(u,\alpha
_{k})\right.  \\
&&-\left. \mathbf{1}_{i,q-k}^{y_{i}}\mathbf{1}_{i,k}^{\alpha _{i}}\int_{%
\mathbf{I}_{i,k}^{y}}du\partial _{1}K(u,y_{1})\ldots \partial
_{1}K(u,y_{q-k})\partial _{1}K(u,\alpha _{1})\ldots \partial _{1}K(u,\alpha
_{k})\right)  \\
&&\hspace*{-0.2in} \hspace*{-0.15in}\left( \mathbf{1}_{i+1,q-k}^{z_{i}}\mathbf{1}_{i+1,k}^{\alpha
_{i}}\int_{\mathbf{I}_{i+1,k}^{z}}^{\frac{i+1}{N}}dv\partial
_{1}K(v,z_{1})\ldots \partial _{1}K(v,z_{q-k})\partial _{1}K(v,\alpha
_{1})\ldots \partial _{1}K(v,\alpha _{k})\right.  \\
&&-\left. \mathbf{1}_{i,q-k}^{y_{i}}\mathbf{1}_{i,k}^{\alpha _{i}}\int_{%
\mathbf{I}_{i,k}^{z}}^{\frac{i}{N}}dv\partial _{1}K(v,z_{1})\ldots \partial
_{1}K(v,z_{q-k})\partial _{1}K(v,\alpha _{1})\ldots \partial _{1}K(v,\alpha
_{k})\right)
\end{eqnarray*}%
where $\mathbf{1}_{i,k}^{x_{j}}$ denotes the indicator function $1_{[0,\frac{%
i}{N}]^{k}}(x_{j})$ with $x$ being $y,z,$ or $\alpha $, and $\mathbf{I}%
_{i,k}^{x}$ denotes the interval $[x_{1}\vee \cdots \vee x_{q-k}\vee \alpha
_{1}\cdots \vee \alpha _{k}$;$i/N]$, with $x$ being $y$ or $z$. By interchanging the order of the integration we get
 \begin{eqnarray}
&&\left\langle f_{i,N}\otimes _{k}f_{i,N})(y_{1},\ldots y_{q-k},z_{1},\ldots ,z_{q-k}\right\rangle \notag\\
&=&d(H,q)^{2}\bigg\{ 1_{[0,\frac{i+1}{N}]^{2q-2k}}(y_{i},z_{i})%
\int_{y_{1}\vee \ldots y_{q-k}}^{\frac{i+1}{N}}du\partial
_{1}K(u,y_{1})\ldots \partial _{1}K(u,y_{q-k})\notag \\
&&\int_{z_{1}\vee \ldots z_{q-k}}^{\frac{i+1}{N}}dv\partial
_{1}K(v,z_{1})\ldots \partial _{1}K(v,z_{q-k})\left( \int_{0}^{u\wedge
v}\partial _{1}K(u,\alpha )\partial _{1}K(v,\alpha )d\alpha \right) ^{k} \notag\\
&&-1_{[0,\frac{i+1}{N}]^{q-k}}(y_{i})1_{[0,\frac{i}{N}]^{q-k}}(z_{i})%
\int_{y_{1}\vee \ldots y_{q-k}}^{\frac{i+1}{N}}du\partial
_{1}K(u,y_{1})\ldots \partial _{1}K(u,y_{q-k}) \notag \\
&&\int_{z_{1}\vee \ldots z_{q-k}}^{\frac{i}{N}}dv\partial
_{1}K(v,z_{1})\ldots \partial _{1}K(v,z_{q-k})\left( \int_{0}^{u\wedge
v}\partial _{1}K(u,\alpha )\partial _{1}K(v,\alpha )d\alpha \right) ^{k}\notag \\
&&-1_{[0,\frac{i}{N}]^{q-k}}(y_{i})1_{[0,\frac{i+1}{N}]^{q-k}}(z_{i})%
\int_{y_{1}\vee \ldots y_{q-k}}^{\frac{i}{N}}du\partial _{1}K(u,y_{1})\ldots
\partial _{1}K(u,y_{q-k}) \notag\\
&&\int_{z_{1}\vee \ldots z_{q-k}}^{\frac{i+1}{N}}dv\partial
_{1}K(v,z_{1})\ldots \partial _{1}K(v,z_{q-k})\left( \int_{0}^{u\wedge
v}\partial _{1}K(u,\alpha )\partial _{1}K
(v,\alpha )d\alpha \right) ^{k}\notag \\
&&+1_{[0,\frac{i}{N}]^{q-k}}(y_{i})1_{[0,\frac{i}{N}]^{q-k}}(z_{i})%
\int_{y_{1}\vee \ldots y_{q-k}}^{\frac{i}{N}}du\partial _{1}K(u,y_{1})\ldots
\partial _{1}K(u,y_{q-k}) \notag \\
&&\int_{z_{1}\vee \ldots z_{q-k}}^{\frac{i}{N}}dv\partial _{1}K(v,z_{1})\ldots \partial _{1}K(v,z_{q-k})\left(
\int_{0}^{u\wedge v}\partial _{1}K(u,\alpha )\partial _{1}K(v,\alpha )d\alpha \right) ^{k}\bigg\}\notag
\end{eqnarray}
and since
$$
\int_{0}^{u\wedge v}\partial _{1}K(u,\alpha )\partial
_{1}K(v,\alpha )d\alpha =a(H^{\prime })|u-v|^{2H^{\prime }-2}
$$
with $a(H^{\prime })=H^{\prime }(2H^{\prime }-1)$, we obtain
\begin{eqnarray}
&&\left\langle f_{i,N}\otimes _{k}f_{i,N}\right\rangle (y_{1},\ldots
y_{q-k},z_{1},\ldots ,z_{q-k})=d(H,q)^{2}a(H^{\prime })^{k}  \notag \\
&&\times \bigg\{1_{[0,\frac{i+1}{N}]^{q-k}}(y_{i})1_{[0,\frac{i+1}{N}%
]^{q-k}}(z_{i})\int_{y_{1}\vee \ldots y_{q-k}}^{\frac{i+1}{N}}du\partial
_{1}K(u,y_{1})\ldots \partial _{1}K(u,y_{q-k})  \notag
\end{eqnarray}%
\begin{eqnarray}
&&\times \int_{z_{1}\vee \ldots z_{q-k}}^{\frac{i+1}{N}}dv\partial
_{1}K(v,z_{1})\ldots \partial _{1}K(v,z_{q-k})|u-v|^{(2H^{\prime }-2)k}
\notag \\
&&-1_{[0,\frac{i+1}{N}]^{q-k}}(y_{i})1_{[0,\frac{i}{N}]^{q-k}}(z_{i})%
\int_{y_{1}\vee \ldots y_{q-k}}^{\frac{i+1}{N}}du\partial
_{1}K(u,y_{1})\ldots \partial _{1}K(u,y_{q-k})  \notag \\
&&\times \int_{z_{1}\vee \ldots z_{q-k}}^{\frac{i}{N}}dv\partial
_{1}K(v,z_{1})\ldots \partial _{1}K(v,z_{q-k})|u-v|^{(2H^{\prime }-2)k}
\notag \\
&&-1_{[0,\frac{i}{N}]^{q-k}}(y_{i})1_{[0,\frac{i+1}{N}]^{q-k}}(z_{i})%
\int_{y_{1}\vee \ldots y_{q-k}}^{\frac{i}{N}}du\partial _{1}K(u,y_{1})\ldots
\partial _{1}K(u,y_{q-k})  \notag \\
&&\times \int_{z_{1}\vee \ldots z_{q-k}}^{\frac{i+1}{N}}dv\partial
_{1}K(v,z_{1})\ldots \partial _{1}K(v,z_{q-k})|u-v|^{(2H^{\prime }-2)k}
\notag \\
&&+1_{[0,\frac{i}{N}]^{q-k}}(y_{i})1_{[0,\frac{i}{N}]^{q-k}}(z_{i})%
\int_{y_{1}\vee \ldots y_{q-k}}^{\frac{i}{N}}du\partial _{1}K(u,y_{1})\ldots
\partial _{1}K(u,y_{q-k})  \notag \\
&&\times \int_{z_{1}\vee \ldots z_{q-k}}^{\frac{i}{N}}dv\partial
_{1}K(v,z_{1})\ldots \partial _{1}K(v,z_{q-k})|u-v|^{(2H^{\prime }-2)k}%
\bigg\}.  \notag \\
&&  \label{fik}
\end{eqnarray}%
As a consequence, we can write
\begin{equation}
V_{N}=T_{2q}+c_{2q-2}T_{2q-2}+\ldots+c_{4}T_{4}+c_{2}T_{2}  \label{VNchaos}
\end{equation}
where
\begin{equation}
c_{2q-2k}:=k!\binom{q}{k}^{2}  \label{comboconst}
\end{equation}
are the combinatorial constants from the product formula for $0\leq k\leq
q-1 $, and%
\begin{equation}
T_{2q-2k}:=N^{2H-1}I_{2q-2k}\left(
\sum_{i=0}^{N-1}f_{i,N}\otimes_{k}f_{i,N}\right) , \label{t2k}
\end{equation}
where the integrands in the last formula above are given explicitly in (\ref%
{fik}). This Wiener chaos decomposition of $V_{N}$ allows us to find $V_{N}$%
's precise order of magnitude via its variance's asymptotics.
\begin{proposition}
\label{PropLim}With
\begin{equation}
c_{1,H}=\frac{4d(H,q)^{4}(H^{\prime}(2H^{\prime}-1))^{2q}}{(4H^{\prime
}-3)(4H^{\prime}-2)[(2H^{'}-2)(q-1)+1]^{2}[(H'-1)(q-1)+1]^{2}},  \label{c1H}
\end{equation}
it holds that%
\begin{equation*}
\lim_{N\rightarrow\infty}\mathbf{E}\left[ c_{1,H}^{-1}N^{(2-2H^{%
\prime})2}c_{2}^{-2}V_{N}^{2}\right] =1.
\end{equation*}
\end{proposition}
\textbf{Proof.} We only need to estimate the $L^{2}$ norm of each term appearing in the chaos decomposition (\ref{VNchaos}) of
$V_{N}$, since these terms are orthogonal in $L^{2}$.
We can write, for $0\leq k\leq q-1$,
\begin{eqnarray*}
\mathbf{E}\left[ T_{2q-2k}^{2}\right]  &=&N^{4H-2}(2q-2k)!\left\Vert \left(
\sum_{i=0}^{N-1}f_{i,N}\otimes _{k}f_{i,N}\right) ^{s}\right\Vert
_{L^{2}([0,1]^{2q-2k})}^{2} \\
&=&N^{4H-2}(2q-2k)!\sum_{i,j=0}^{N-1}\langle f_{i,N}\tilde{\otimes}%
_{k}f_{i,N},f_{j,N}\tilde{\otimes}_{k}f_{j,N}\rangle _{L^{2}([0,1]^{2q-2k})}
\end{eqnarray*}%
where $(g)^{s}=\tilde{g}$ and $f_{i,N}\tilde{\otimes}_{k}f_{i,N}$ denotes
the symmetrization of the function $f_{i,N}\otimes _{k}f_{i,N}.$ We will
consider first the term $T_{2}$ obtained for $k=q-1$. In this case, the
kernel $\sum_{i=0}^{N-1}f_{i,N}\otimes _{q-1}f_{i,N}$ is symmetric and we
can avoid its symmetrization. Therefore
\begin{eqnarray*}
\mathbf{E}\left[ T_{2}^{2}\right]  &=&2!\;N^{4H-2}\Vert
\sum_{i=0}^{N-1}f_{i,N}\otimes _{q-1}f_{i,N}\Vert _{L^{2}([0,1]^{2})}^{2} \\
&=&2!\;N^{4H-2}\sum_{i,j=0}^{N-1}\langle f_{i,N}\otimes
_{q-1}f_{i,N},f_{j,N}\otimes _{q-1}f_{j,N}\rangle _{L^{2}([0,1]^{2})}.
\end{eqnarray*}%
We compute now the scalar product in the above expression. By using Fubini
theorem, we end up with the following easier expression
\begin{eqnarray*}
&&\langle f_{i,N}\otimes _{q-1}f_{i,N},f_{j,N}\otimes _{q-1}f_{j,N}\rangle
_{L^{2}([0,1]^{2})}=a(H^{\prime
})^{2q}d(H,q)^{4}\int_{I_{i}}\int_{I_{i}}\int_{I_{j}}\int_{I_{j}} \\
&&\times |u-v|^{(2H^{\prime }-2)(q-1)}|u^{\prime }-v^{\prime }|^{(2H^{\prime
}-2)(q-1)}|u-u^{\prime }|^{2H^{\prime }-2}|v-v^{\prime }|^{2H^{\prime
}-2}dv^{\prime }du^{\prime }dvdu
\end{eqnarray*}%
Using the change of variables $y=(u-\frac{i}{N})N$ and similarly for the
other variables, we now obtain
\begin{eqnarray*}
&&\mathbf{E}\left[ T_{2}^{2}\right] =2d(H,q)^{4}(H^{\prime }(2H^{\prime
}-1))^{2q}N^{4H-2}N^{-4}N^{-(2H^{\prime }-2)2q} \\
&&\times
\sum_{i,j=0}^{N-1}\int_{0}^{1}\int_{0}^{1}\int_{0}^{1}\int_{0}^{1}dydzdy^{%
\prime }dz^{\prime }|y-z|^{(2H^{\prime }-2)(q-1)}|y^{\prime }-z^{\prime
}|^{(2H^{\prime }-2)(q-1)} \\
&&\times |y-y^{\prime }+i-j|^{(2H^{\prime }-2)}|z-z^{\prime
}+i-j|^{(2H^{\prime }-2)}.
\end{eqnarray*}%
This can be viewed as the sum of a diagonal part ($i=j$) and a non-diagonal
part ($i\neq j$), where the non-diagonal part is dominant, as the reader
will readily check. Therefore, the behavior of $\mathbf{E}\left[ T_{2}^{2}%
\right] $ will be given by
\begin{eqnarray*}
&&\mathbf{E}\left[ T_{2}^{\prime 2}\right]  \\
&:=&2!d(H,q)^{4}(H^{\prime }(2H^{\prime
}-1))^{2q}N^{-2}2\sum_{i>j}\int_{0}^{1}\int_{0}^{1}\int_{0}^{1}%
\int_{0}^{1}dydzdy^{\prime }dz^{\prime } \\
&&\times (|y-z||y^{\prime }-z^{\prime }|)^{(2H^{\prime
}-2)(q-1)}(|y-y^{\prime }+i-j||z-z^{\prime }+i-j|)^{(2H^{\prime }-2)} \\
&=&2!d(H,q)^{4}(H^{\prime }(2H^{\prime
}-1))^{2q}N^{-2}2\sum_{i=0}^{N-2}\sum_{\ell
=2}^{N-i}\int_{0}^{1}\int_{0}^{1}\int_{0}^{1}\int_{0}^{1}dydzdy^{\prime
}dz^{\prime } \\
&&\times (|y-z||y^{\prime }-z^{\prime }|)^{(2H^{\prime
}-2)(q-1)}(|y-y^{\prime }+\ell -1||z-z^{\prime }+\ell -1|)^{(2H^{\prime }-2)}
\\
&=&2!d(H,q)^{4}(H^{\prime }(2H^{\prime }-1))^{2q}N^{-2}2\sum_{\ell
=2}^{N}(N-\ell
+1)\int_{0}^{1}\int_{0}^{1}\int_{0}^{1}\int_{0}^{1}dydzdy^{\prime
}dz^{\prime } \\
&&\times (|y-z||y^{\prime }-z^{\prime }|)^{(2H^{\prime
}-2)(q-1)}(|y-y^{\prime }+\ell -1||z-z^{\prime }+\ell -1|)^{(2H^{\prime }-2)}.
\end{eqnarray*}%
As in \cite{TV} note that
\begin{eqnarray*}
&&\frac{1}{N^{2}}\sum_{\ell =2}^{N}(N-\ell +1)|y-y^{\prime }+\ell
-1|^{(2H^{\prime }-2)}|z-z^{\prime }+\ell -1|^{(2H^{\prime }-2)} \\
&=&N^{2(2H^{\prime }-2)}\frac{1}{N}\sum_{\ell =2}^{N}(1-\frac{\ell -1}{N})|%
\frac{y-y^{\prime }}{N}+\frac{\ell -1}{N}|^{2H^{\prime }-2}|\frac{%
z-z^{\prime }}{N}+\frac{\ell -1}{N}|^{2H^{\prime }-2}.
\end{eqnarray*}%
Using a Riemann sum approximation argument we conclude that
\begin{equation*}
\mathbf{E}\left[ T_{2}^{^{\prime }2}\right] \sim \frac{4d(H,q)^{4}(H^{\prime
}(2H^{\prime }-1))^{2q}\;\times \;N^{2(2H^{\prime }-2)}}{(4H^{\prime
}-3)(4H^{\prime }-2)[((2H^{^{\prime }}-2)(q-1)+1)]^{2}[(H^{\prime
}-1)(q-1)+1]^{2}}.
\end{equation*}%
Therefore, it follows that
\begin{equation}
\mathbf{E}\left[ c_{1,H}^{-1}N^{2(2-2H^{\prime })}T_{2}^{2}\right]
\rightarrow _{N\rightarrow \infty }1,  \label{meanT2}
\end{equation}%
with $c_{1,H}$ as in (\ref{c1H}). \\
Let us study now the term $T_{4},\ldots
,T_{2q}$ given by (\ref{t2k}). Here the function \\$\sum_{i=0}^{N-1}f_{i,N}%
\otimes _{k}f_{i,N}$ is no longer symmetric but we will show that the
behavior of its $L^{2}$ norm is dominated by $\mathbf{E}\left[ T_{2}^{2}%
\right] $. Since for any square integrable function $g$ one has $\Vert
\tilde{g}\Vert _{L^{2}}\leq \Vert g\Vert _{L^{2}}$, we have for $k=0,\ldots
,q-2$
\begin{eqnarray}
&& \frac{1}{(2q-2k)!}\mathbf{E}\left[ T_{2q-2k}^{2}\right]  =N^{4H-2}\Vert
\sum_{i=0}^{N-1}f_{i,N}\tilde{\otimes}_{k}f_{i,N}\Vert
_{L^{2}([0,1]^{2q-2k})}^{2}  \notag \\
&\leq &N^{4H-2}\Vert \sum_{i=0}^{N-1}f_{i,N}\otimes _{k}f_{i,N}\Vert
_{L^{2}([0,1]^{2q-2k})}^{2}  \notag \\
&=&N^{4H-2}\sum_{i,j=0}^{N-1}\langle f_{i,N}\otimes
_{k}f_{i,N},f_{j,N}\otimes _{k}f_{j,N}\rangle _{L^{2}([0,1]^{2q-2k})}
\label{forpolarization}
\end{eqnarray}%
and proceeding as above, with $e_{H,q,k}:=(2q-2k)!(H^{\prime }(2H^{\prime
}-1))^{2q}d(H,q)^{4}$ we can write
\begin{eqnarray*}
&&\mathbf{E}\left[ T_{2q-2k}^{2}\right] \leq
e_{H,q,k}N^{4H-2}\sum_{i,j=0}^{N-1}\int_{I_{i}}\int_{I_{i}}dy_{1}dz_{1}%
\int_{I_{j}}\int_{I_{j}}dy_{1}^{\prime }dz_{1}^{\prime } \\
&&\times |y_{1}-z_{1}|^{(2H^{\prime }-2)k}|y_{1}^{\prime }-z_{1}^{\prime
}|^{(2H^{\prime }-2)k}|y_{1}-y_{1}^{\prime }|^{(2H^{\prime
}-2)(q-k)}|z_{1}-z_{1}^{\prime }|^{(2H^{\prime }-2)(q-k)}
\end{eqnarray*}%
and using a change of variables as before,%
\begin{eqnarray}
&&\mathbf{E}\left[ T_{2q-2k}^{2}\right]   \notag \\
&\leq &e_{H,q,k}N^{4H-2}N^{-4}N^{-(2H^{\prime
}-2)2q}\sum_{i,j=0}^{N-1}\int_{[0,1]^{4}}dydzdy^{\prime }dz^{\prime
}(|y-z||y^{\prime }-z^{\prime }|)^{(2H^{\prime }-2)k}  \notag \\
&&\times |y-y^{\prime }+i-j|^{(2H^{\prime }-2)(q-k)}|z-z^{\prime
}+i-j|^{(2H^{\prime }-2)(q-k)}  \notag \\
&=&\frac{(2q-2k)!d(H,q)^{4}}{a(H^{\prime })^{-2q}}\frac{N^{(2H^{\prime
}-2)(2q-2k)}}{N^{2}}\sum_{\ell =2}^{N}\left( 1-\left( \frac{\ell -1}{N}%
\right) \right) \int_{[0,1]^{4}}dydzdy^{\prime }dz^{\prime }  \notag \\
&&\times (|y-z||y^{\prime }-z^{\prime }|)^{(2H^{\prime }-2)k}\left(
\left\vert \frac{y-y^{\prime }}{N}+\frac{\ell -1}{N}\right\vert \left\vert
\frac{z-z^{\prime }}{N}+\frac{\ell -1}{N}\right\vert \right) ^{(2H^{\prime
}-2)(q-k)}.  \label{T2q}
\end{eqnarray}%
Since off a diagonal term (again of lower order), the terms $\frac{%
z-z^{\prime }}{N}$ are dominated by $\frac{\ell }{N}$ for large $l,N$ it
follows that, for $1\leq k\leq q-1$%
\begin{equation}
\mathbf{E}\left[ c_{q-k,H}^{-1}N^{(2-2H^{\prime })(2q-2k)}T_{2q-2k}^{2}%
\right] ={O}(1)  \label{meanT2k}
\end{equation}%
when $N\rightarrow \infty $, with
\begin{equation}
c_{q-k,H}=2(\int_{0}^{1}(1-x)x^{(2H^{\prime }-2)(2q-2k)}dx)a(H^{\prime
})^{-2}(2q-2k)!d(H,q)^{2}a(H^{\prime })^{2q}.  \label{c1Hq}
\end{equation}%
It is obvious that the dominant term in the decomposition of $V_{N}$ is the
term in the chaos of order 2. [The case $k=0$ is in the same situation for $%
H>\frac{3}{4}$ and for $H\in (\frac{1}{2},\frac{3}{4})$ the term $T_{2q}$
obtained for $k=0$ has to be renormalized by $N$, (as in proofs in the
remainder of the paper); in any case it is dominated by the term $T_{2}$]. More specifically we have for any $k\leq q-2$,
\begin{equation}
\mathbf{E}\left[ N^{2(2-2H^{\prime})}T_{2q-2k}^{2}\right] =O\left(
N^{-2\left( 2-2H^{\prime}\right) 2\left( q-k-1\right) }\right) .
\label{lower}
\end{equation}
Combining this with the orthogonality of chaos integrals, we immediately get
that, up to terms that tend to $0$, $N^{2-2H^{\prime}}V_{N}$ and $%
N^{2-2H^{\prime}}T_{2}$ have the same norm in $L^{2}\left( \Omega\right) $.
This finishes the proof of the proposition.\hfill\vrule width.25cm
height.25cm depth0cm\smallskip
\vskip0.5cm
Summarizing the spirit of the above proof, to understand the behavior of the
renormalized sequence $V_{N}$ it suffices to study the limit of the term%
\begin{equation}
I_{2}\left( N^{2H-1}N^{(2-2H^{\prime})}\sum_{i=0}^{N-1}f_{i,N}\otimes
_{q-1}f_{i,N}\right)  \label{i2}
\end{equation}
with
\begin{eqnarray}
&&(f_{i,N}\otimes _{q-1}f_{i,N})(y,z)    \label{fourbrothers} \\
&=&d(H,q)^{2}a(H^{\prime })^{q-1}  \notag \\
&&\left( 1_{[0,\frac{i}{N}]}(y\vee z)\int_{I_{i}}\int_{I_{i}}dvdu\partial
_{1}K(u,y)\partial _{1}K(v,z)|u-v|^{(2H^{\prime }-2)(q-1)}\right.   \notag \\
&&\left. +1_{[0,\frac{i}{N}]}(y)1_{I_{i}}(z)\int_{I_{i}}\int_{z}^{\frac{i+1}{%
N}}dvdu\partial _{1}K(u,y)\partial _{1}K(v,z)|u-v|^{(2H^{\prime
}-2)(q-1)}\right.   \notag \\
&&\left. +1_{I_{i}}(y)1_{[0,\frac{i}{N}]}(z)\int_{y}^{\frac{i+1}{N}%
}\int_{I_{i}}dvdu\partial _{1}K(u,y)\partial _{1}K(v,z)|u-v|^{(2H^{\prime
}-2)(q-1)}\right.   \notag \\
&&\left. +1_{I_{i}}(y)1_{I_{i}}(z)\int_{y}^{\frac{i+1}{N}}\int_{z}^{\frac{i+1%
}{N}}dvdu\partial _{1}K(u,y)\partial _{1}K(v,z)|u-v|^{(2H^{\prime }-2)(q-1)}\right) .\notag
\end{eqnarray}%
We will see in the proof of the next theorem that, of the contribution of
the four terms on the right-hand side of (\ref{fourbrothers}), only the
first one does not tend to $0$ in $L^{2}\left( \Omega \right) $. Hence the
following notation will be useful: $f_{2}^{N}$ will denote the integrand of
the contribution to (\ref{i2}) corresponding to that first term, and $r_{2}$
will be the remainder of the integrand in (\ref{i2}). In other words,
\begin{equation}
f_{2}^{N}+r_{2}=N^{2H-1}N^{(2-2H^{\prime })}\sum_{i=0}^{N-1}f_{i,N}\otimes
_{q-1}f_{i,N}  \label{r2d2c3po}
\end{equation}%
and%
\begin{eqnarray}
f_{2}^{N}(y,z) &:&=N^{2H-1}N^{(2-2H^{\prime })}d(H,q)^{2}a(H^{\prime })^{q-1}
\notag \\
&&\hspace*{-1in}\sum_{i=0}^{N-1}1_{[0,\frac{i}{N}]}(y\vee
z)\int_{I_{i}}\int_{I_{i}}dvdu\partial _{1}K(u,y)\partial
_{1}K(v,z)|u-v|^{(2H^{\prime }-2)(q-1)}.  \label{f2N}
\end{eqnarray}
\begin{theorem}
\label{ThmConv}The sequence given by (\ref{i2}) converges in $L^{2}\left(
\Omega\right) $ as $N\rightarrow\infty$ to the constant $c_{1,H}^{1/2}$ times a standard Rosenblatt random variable $Z_{1}^{(2,2H^{\prime
}-1)}$ with selfsimilarity parameter $2H^{\prime}-1$ and $H^{\prime}$ is given
by (\ref{H'}). Consequently, we also have that $c_{1,H}^{-1/2}N^{(2-2H^{%
\prime})}c_{2}^{-1}V_{N}$ converges in $L^{2}\left( \Omega\right) $ as $%
N\rightarrow\infty$ to the same Rosenblatt random variable.
\end{theorem}
\textbf{Proof: } The first statement of the theorem is that $ N^{2-2H^{\prime}}T_{2}$ converges to $$c_{1,H}^{1/2}Z_{1}^{(2,2H^{\prime}-1)}$$ in $L^{2}(\Omega)$. From (\ref{i2}) it follows that $T_{2}$ is a second-chaos random
variable, with kernel $$ N^{2H-1}\sum_{i=0}^{N-1}(f_{i,N}\otimes _{q-1}f_{i,N})$$ (see expression in ( \ref{fourbrothers})),
so we only need to prove this kernel converges in $ L^{2}\left( [0,1]^{2}\right) $.
The first observation is  that $r_{2}(y,z)$ defined in (\ref{r2d2c3po})
converges to zero in $L^{2}([0,1]^{2})$ as $N\rightarrow \infty $. The
crucial fact is that the intervals $I_{i}$ which are disjoints, appear in
the expression of this term and this implies that the non-diagonal terms
vanish when we take the square norm of the sum; in fact it can easily be seen that the norm
in $L^{2}$ of $r_{2}$ corresponds to the diagonal part in the evaluation in $%
ET_{2}^{2}$ which is clearly dominated by the non-diagonal part, so this
result comes as no surprise. The proof follows the lines of \cite{TV}.
This shows $N^{(2-2H^{\prime })}T_{2}$ is the sum of $I_{2}\left(
f_{2}^{N}\right) $ and a term which goes to $0$ in $L^{2}\left( \Omega
\right) $. Our next step is thus simply to calculate the limit in $%
L^{2}\left( \Omega \right) $, if any, of $I_{2}(f_{2}^{N})$ where $f_{2}^{N}$
has been defined in (\ref{f2N}). By the isometry property (\ref{isom}),
limits of second-chaos r.v.'s in $L^{2}\left( \Omega \right) $ are
equivalent to limits of their symmetric kernels in $L^{2}\left(
[0,1]^{2}\right) $. Note that $f_{2}^{N}$ is symmetric. Therefore, it is
sufficient to prove that $f_{2}^{N}$ converges to the kernel of the
Rosenblatt process at time 1. We have by definition%
\begin{eqnarray*}
f_{2}^{N}(y,z) &=&(H^{\prime }(2H^{\prime
}-1))^{(q-1)}d(H,q)^{2}N^{2H-1}N^{2-2H^{\prime }} \\
&&\times \sum_{i=0}^{N-1}\int_{I_{i}}\int_{I_{i}}|u-v|^{(2H^{\prime
}-2)(q-1)}\partial _{1}K^{H^{\prime }}(u,y)\partial _{1}K^{H^{\prime }}(v,z).
\end{eqnarray*}
Thus for every $y,z$,
\begin{eqnarray*}
&&f_{2}^{N}(y,z)=d(H,q)^{2}(H^{\prime }(2H^{\prime
}-1))^{(q-1)}N^{2H-1}N^{2-2H^{\prime }} \\
&&\times \sum_{i=0}^{N-1}\int_{I_{i}}\int_{I_{i}}|u-v|^{(2H^{\prime
}-2)(q-1)}\partial _{1}K^{H^{\prime }}(u,y)\partial _{1}K^{H^{\prime
}}(v,z)dudv \\
&=&d(H,q)^{2}(H^{\prime }(2H^{\prime }-1))^{(q-1)}N^{2H-1}N^{2-2H^{\prime
}}\sum_{i=0}^{N-1}\int_{I_{i}}\int_{I_{i}}|u-v|^{(2H^{\prime }-2)(q-1)} \\
&&\times \left( \partial _{1}K^{H^{\prime }}(u,y)\partial _{1}K^{H^{\prime
}}(v,z)-\partial _{1}K^{H^{\prime }}(i/N,z)\partial _{1}K^{H^{\prime
}}(i/N,z)\right) dudv \\
&&+d(H,q)^{2}(H^{\prime }(2H^{\prime }-1))^{(q-1)}N^{2H-1}N^{2-2H^{\prime
}}\sum_{i=0}\int_{I_{i}}\int_{I_{i}}|u-v|^{(2H^{\prime }-2)(q-1)} \\
&&\times \partial _{1}K^{H^{\prime }}(i/N,y)\partial _{1}K^{H^{\prime
}}(i/N,z)dudv \\
&=:&A_{1}^{N}(y,z)+A_{2}^{N}(y,z).
\end{eqnarray*}%
As in \cite{TV}, one can show that 
$\mathbf{E}\left[ \left\Vert A_{1}^{N}\right\Vert _{L^{2}([0,1]^{2})}^{2}%
\right] \rightarrow 0\text{ as }N\rightarrow \infty $. 
Regarding the second term $A_{2}^{N}(y,z)$, the summand is zero if $%
i/N<y\vee z$, therefore we get that $f_{2}^{N}$ is equivalent to
\begin{eqnarray*}
&&N^{2H-1}N^{2-2H^{\prime }}d(H,q)^{2}(H^{\prime }(2H^{\prime
}-1))^{(q-1)}\sum_{i=0}^{N-1}\int_{I_{i}}\int_{I_{i}}|u-v|^{(2H^{\prime
}-2)(q-1)} \\
&&\times \partial _{1}K^{H^{\prime }}(i/N,y)\partial _{1}K^{H^{\prime
}}(i/N,z)dudv
\end{eqnarray*}%
\begin{eqnarray*}
&=&(H^{\prime }(2H^{\prime }-1))^{(q-1)}d(H,q)^{2}N^{2H-1}N^{2-2H^{\prime }}
\\
&&\times \sum_{i=0}^{N-1}\partial _{1}K^{H^{\prime }}(i/N,y)\partial
_{1}K^{H^{\prime }}(i/N,z)1_{y\vee
z<i/N}\int_{I_{i}}\int_{I_{i}}|u-v|^{(2H^{\prime }-2)(q-1)}dudv \\
&=&(H^{\prime }(2H^{\prime }-1))^{(q-1)}[((2H^{\prime
}-2)(q-1)+1)((H^{\prime }-1)(q-1)+1)]^{-1} \\
&&\times \frac{N^{2H-1}N^{(2-2H^{\prime })q}}{N^{2}}\sum_{i=0}^{N-1}\partial
_{1}K^{H^{\prime }}(i/N,y)\partial _{1}K^{H^{\prime }}(i/N,z)1_{y\vee z<i/N}
\\
&=&\frac{d(H,q)^{2}}{d(2H^{\prime }-2,2)}\frac{(H^{\prime }(2H^{\prime
}-1))^{(q-1)}}{((2H^{\prime }-2)(q-1)+1)((H^{\prime }-1)(q-1)+1)} \\
&&\times d(2H^{\prime }-2,2)N^{-1}\sum_{i=0}^{N-1}\partial _{1}K^{H^{\prime
}}(i/N,y)\partial _{1}K^{H^{\prime }}(i/N,z)1_{y\vee z<i/N}.
\end{eqnarray*}%
The sequence $d(2H^{\prime }-2,2)N^{-1}\sum_{i=0}^{N-1}\partial
_{1}K^{H^{\prime }}(i/N,y)\partial _{1}K^{H^{\prime }}(i/N,z)1_{y\vee z<i/N}$
is a Riemann sum that converges pointwise on $[0,1]^{2}$ to the kernel of
the Rosenblatt process $Z^{2H^{\prime }-1,2}$ at time $1$. To obtain the
convergence in $L^{2}\left( [0,1]^{2}\right) $ we will apply the dominated
convergence theorem. Indeed,
\begin{eqnarray*}
&&\int_{0}^{1}\int_{0}^{1}\left\vert \frac{1}{N}\sum_{i=0}^{N-1}\partial
_{1}K^{H^{\prime }}(i/N,y)\partial _{1}K^{H^{\prime }}(i/N,z)1_{y\vee
z<i/N}\right\vert ^{2}dydz \\
&=&\frac{1}{N^{2}}\sum_{i,j=0}^{N-1}\left\vert \int_{0}^{1}\partial
_{1}K^{H^{\prime }}(i/N,y)\partial _{1}K^{H^{\prime }}(j/N,y)1_{y<(i\wedge
j)/N}dy\right\vert ^{2} \\
&\leq &\frac{1}{N^{2}}\sum_{i,j=0}^{N-1}\left\vert \mathbf{E}\left[ \Delta
Z_{i/N}\Delta Z_{j/N}\right] \right\vert ^{2},
\end{eqnarray*}%
where $\Delta Z_{i/N}$ is the difference $Z\left( \frac{i}{N}\right)
-Z\left( \frac{i-1}{N}\right) $ for a Rosenblatt process $Z$. We now show
that the above sum is always $\ll $ $N^{2}$, which proves that the last
expression, with the $N^{-2}$ factor, is bounded. In fact for $%
H_{1}=2H^{\prime }-1$
\begin{eqnarray*}
&&\hspace*{-0.35in}\sum_{i,j=0}^{N-1}\left\vert \mathbf{E}\left[ \Delta Z_{i/N}\Delta Z_{j/N}%
\right] \right\vert ^{2}=\sum_{i,j=0}^{N-1}\left\vert \left\vert \frac{i-j+1%
}{N}\right\vert ^{2H_{1}}+\left\vert \frac{i-j-1}{N}\right\vert
^{2H_{1}}-2\left\vert \frac{i-j}{N}\right\vert ^{2H_{1}}\right\vert ^{2} \\
&=&\frac{N^{-4H_{1}}}{4}\sum_{i,j=0}^{N-1}\left\vert \left\vert
i-j+1\right\vert ^{2H_{1}}+\left\vert i-j-1\right\vert ^{2H_{1}}-2\left\vert
i-j\right\vert ^{2H_{1}}\right\vert ^{2} \\
&\leq &\frac{N^{-4H_{1}}}{4}2N\sum_{\ell =-N+1}^{N-1}\left\vert \left\vert
\ell +1\right\vert ^{2H_{1}}+\left\vert \ell -1\right\vert
^{2H_{1}}-2\left\vert \ell \right\vert ^{2H_{1}}\right\vert ^{2}.
\end{eqnarray*}%
The function $g(\ell )=|\ell +1|^{2H_{1}}+|\ell -1|^{2H_{1}}-2|\ell
|^{2H_{1}}$ behaves like $H_{1}(2H_{1}-1)|\ell |^{2H_{1}-2}$ for large $\ell $. We need to
separate the cases of convergence and divergence of the series $%
\sum_{-\infty }^{\infty }|g(\ell )|^{2}$. It is divergent as soon as $%
H_{1}\geq 3/4$, or equivalently $H^{\prime }\geq 7/8$, in which case we get
for some constant $c$ not dependent on $N$,
\begin{equation*}
\sum_{i,j=0}^{N-1}\left\vert \mathbf{E}\left[ \Delta Z_{i/N}\Delta Z_{j/N}%
\right] \right\vert ^{2}\leq cN^{-4H_{1}+1+4H_{1}-3}=cN^{-2}\ll N^{2}.
\end{equation*}%
The series $\sum_{-\infty }^{\infty }|g(\ell )|^{2}$ is convergent if $%
H^{\prime }<7/8$, in which case we get%
\begin{equation*}
\sum_{i,j=0}^{N-1}\left\vert \mathbf{E}\left[ \Delta Z_{i/N}\Delta Z_{j/N}%
\right] \right\vert ^{2}\leq cN^{-4H_{1}+1}.
\end{equation*}%
For this to be $\ll N^{2}$, we simply need $-4H_{1}+1<2$, i.e. $H^{\prime
}>5/8$. However, since $q\geq 2$ and $H>1/2$ we always have $H^{\prime }>3/4$%
. Therefore in all cases, the sequence $A_{2}^{N}(y,z)$ is bounded in $%
L^{2}\left( [0,1]^{2}\right)$ and in this way we obtain the $L^{2}$
convergence to the kernel of a Rosenblatt process of order 1. The first
statement of the theorem is proved.
In order to show that $c_{1,H}^{-1/2}N^{(2-2H^{^{\prime }})}c_{2}^{-1}V_{N}$
converges in $L^{2}(\Omega )$ to the same Rosenblatt random variable as the
normalized version of the quantity in (\ref{i2}), it is sufficient to show
that, after normalization by $N^{2-2H^{\prime }}$, each of the remaining
terms of in the chaos expansion (\ref{VNchaos}) of $V_{N}$, converge to zero
in $L^{2}(\Omega )$, i.e. that $N^{(2-2H^{^{\prime }})}T_{2q-2k}$ converge
to zero in $L^{2}(\Omega )$, for all $1\leq k<q-1$. From (\ref{lower}) we
have%
\begin{equation*}
\mathbf{E}\left[ N^{2(2-2H^{\prime })}T_{2q-2k}^{2}\right] =O\left(
N^{-2\left( 2-2H^{\prime }\right) 2\left( q-k-1\right) }\right)
\end{equation*}%
which is all that is needed, concluding the proof of the theorem. \hfill %
\vrule width.25cm height.25cm depth0cm

\smallskip
\section{Reproduction property for the Hermite process\label{ReprodSec}}
We now study the limits of the other terms in the chaos expansion (\ref%
{VNchaos}) of $V_{N}$. We will consider first the convergence of the term of
greatest order $T_{2q}$ in this expansion. The behavior of $T_{2q}$ is
interesting because it differs from the behavior of the all other terms: its
suitable normalization possesses a Gaussian limit if $H\in(1/2,3/4]$.
Therefore it inherits, in some sense, the properties of the quadratic
variations for the fractional Brownian motion (results in \cite{TV}).
We have
\begin{equation*}
T_{2q}=N^{2H-1}I_{2q}\left( \sum_{i=0}^{N-1}f_{i,N}\otimes f_{i,N}\right)
\end{equation*}%
and $ \mathbf{E}\left[ T_{2q}^{2}\right] =N^{4H-2}(2q)!\sum_{i,j=0}^{N-1}\langle
f_{i,N}\tilde{\otimes}f_{i,N},f_{j,N}\tilde{\otimes}f_{j,N}\rangle _{L^{2}[0,1]^{2q}}.$ We will use the following
combinatorial formula: if $f,g$ are two symmetric functions in $L^{2}([0,1]^{q})$,
\begin{eqnarray*}
&&(2q)!\langle f\tilde{\otimes}f,g\tilde{\otimes}g\rangle
_{L^{2}([0,1]^{2q})} \\
&=&(q!)^{2}\langle f\otimes f,g\otimes g\rangle _{L^{2}([0,1]^{2q})}+\sum_{k=1}^{q-1}\binom{q}{k}^{2}(q!)^{2}\langle f\otimes
_{k}g, g\otimes _{k}f\rangle _{L^{2}([0,1]^{2q-2k})}
\end{eqnarray*}%
to obtain
\begin{eqnarray*}
\mathbf{E}\left[ T_{2q}^{2}\right]
&=&N^{4H-2}(q!)^{2}\sum_{i,j=0}^{N-1}\langle f_{i,N},f_{j,N}\rangle
_{L^{2}([0,1]^{q})}^{2} \\
&&+N^{4H-2}(q!)^{2}\sum_{i,j=0}^{N-1}\sum_{k=1}^{q-1}\binom{q}{k}^{2}\langle f_{i,N}\otimes _{k}f_{j,N}, f_{j,N}\otimes
_{k}f_{i,N} \rangle _{L^{2}([0,1]^{2q-2k})}.
\end{eqnarray*}%
First note that $
\mathbf{E}\left[ \left( Z_{\frac{i+1}{N}}^{(q,H)}-Z_{\frac{i}{N}%
}^{(q,H)}\right) \left( Z_{\frac{j+1}{N}}^{(q,H)}-Z_{\frac{j}{N}%
}^{(q,H)}\right) \right]   =\mathbf{E}\left[ I_{q}(f_{i,N})I_{q}(f_{j,N})\right] =q!\langle f_{i,N},f_{j,N}\rangle
_{L^{2}([0,1]^{q})}
$
and so the covariance structure of $Z^{\left( q,H\right) }$ implies%
\begin{eqnarray}
&&(q!)^{2}\sum_{i,j=0}^{N-1}\langle f_{i,N},f_{j,N}\rangle
_{L^{2}([0,1]^{q})}^{2}  \notag \\
&=&\frac{1}{4}\sum_{i,j=0}^{N-1}\left[ \left( \frac{i-j+1}{N}\right)
^{2H}+\left( \frac{i-j-1}{N}\right) ^{2H}-2\left( \frac{i-j}{N}\right) ^{2H}%
\right] ^{2}.  \label{usefBmcov}
\end{eqnarray}%
Secondly, the square norm of the contraction $f_{i,N}\otimes _{k}f_{j,N} $ has been computed before (actually its expression
is obtained in the lines from formula (\ref{forpolarization}) to formula (\ref{T2q})). By a simple
polarization, we obtain
\begin{eqnarray*}
&&N^{2H-1}\sum_{i,j=0}^{N-1}\langle f_{i,N}\otimes
_{k}f_{j,N},f_{j,N}\otimes _{k}f_{i,N}\rangle _{L^{2}([0,1]^{2q-2k})} \\
&&=d(H,q)^{4}a(H^{\prime })^{2q}N^{4H-2}N^{(2H^{\prime
}-2)2q}\sum_{i,j=0}^{N-1}\int_{0}^{1}\int_{0}^{1}\int_{0}^{1}%
\int_{0}^{1}dydzdy^{\prime }dz^{\prime } \\
&&\hspace*{-.25in}(|y-z+i-j||y^{\prime }-z^{\prime }+i-j|)^{(2H^{\prime }-2)k}|(y-y^{\prime
}+i-j||z-z^{\prime }+i-j|)^{(2H^{\prime }-2)(q-k)}
\end{eqnarray*}%
and as in the proof of Proposition \ref{PropLim} we can find that this term
has to be renormalized by, if $H>\frac{3}{4}$
\begin{equation*}
b_{H,k}N^{(2-2H^{\prime })2q}=b_{H,k}N^{4-4H}
\end{equation*}%
where $b_{1,H,k}=(q!)^{2}(C_{k}^{q})^{2}2d(H,q)^{4}a(H^{\prime
})^{2q}\int_{0}^{1}(1-x)x^{4H-4}dx$. If $H\in (\frac{1}{2},\frac{3}{4})$,
then, the same quantity will be renormalized by $b_{2,H,k}N$ where
\begin{equation*}
b_{2,H,k}=(q!)^{2}(C_{k}^{q})^{2}d(H,q)^{4}a(H^{\prime
})^{2q}\sum_{k=1}^{\infty }\left( 2k^{2H}-(k+1)^{2H}-(k-1)^{2H}\right) ^{2}
\end{equation*}
while for $H=\frac{3}{4}$ the renormalization is of order $b_{3,H,k}N(\log
N)^{-1}$ with $b_{3,H,k}=(q!)^{2}(C_{k}^{q})^{2}2d(H,q)^{4}a(H^{\prime
})^{2q}(1/2)$. As a consequence, we find a sum whose behavior is well-known
(it is the same as the mean square of the quadratic variations of the
fractional Brownian motion, see e.g. \cite{TV}) and we get, for large $N$
\begin{eqnarray*}
\mathbf{E}\left[ T_{2q}^{2}\right]  &\sim &\frac{1}{N}x_{1,H},\hskip0.1cm%
\mbox{
if }H\in (\frac{1}{2},\frac{3}{4});\mathbf{E}\left[ T_{2q}^{2}\right] \sim
N^{4H-4}x_{2,H},\hskip0.1cm\mbox{ if }H\in (\frac{3}{4},1) \\
\mathbf{E}\left[ T_{2q}^{2}\right]  &\sim &\frac{\log N}{N}x_{3,H},\hskip%
0.1cm\mbox{ if }H=\frac{3}{4}
\end{eqnarray*}%
where $x_{1,H}=\left( \sum_{l=1}^{q-1}b_{2,H,l}+1+(1/2)\sum_{k=1}^{\infty
}\left( 2k^{2H}-(k+1)^{2H}-(k-1)^{2H}\right) ^{2}\right) $, $%
x_{2,H}=(\sum_{l=1}^{q-1}b_{1,H,l}+H^{2}(2H-1)/(4H-3))$ and $%
x_{3,H}=(\sum_{l=1}^{q-1}b_{3,H,l}+9/16)$.

\vskip0.5cm

\begin{remark}
The fact that the normalizing factor for $T_{2q}$ when $H<3/4$ is $N^{-1/2}$ (in particular does not depend on $H$) is a
tell-tale sign that normal convergence may occur. It is possible to show that the term
in the chaos of order $2q$ converges, after its renormalization, to a Gaussian law. More precisely: a) suppose that $H\in(\frac{1}{2},\frac{3}{4})$ and let $%
F_{N}:=x_{1,H}^{-1/2}\sqrt{N}T_{2q}$; then, as $N\rightarrow\infty$, the
sequence $F_{N}$ convergence to the standard normal distribution $N(0,1)$; b)  suppose $H=\frac{3}{4}$ and set $G_{N}:=\sqrt{\frac{N}{\log N}}%
x_{3,H}^{-1/2}T_{2q}$; then, as $N\rightarrow\infty$, the sequence $G_{N}$ convergence to the standard normal distribution
$N(0,1)$.  We refer to the extended version of our paper on arxiv for the basic ideas of the proof.
\end{remark}

\vskip0.5cm
It is possible to give the limits of the terms $T_{2q-2}$ to $T_{2}$
appearing in the decomposition of the statistics $V_{N}$. All these
renormalized terms will converge to Hermite random variables of the same
order as their indices (we have already proved this property in detail for $%
T_{2}$ in the previous section). This is a kind of \textquotedblleft
reproduction\textquotedblright\ of the Hermite processes through their
variations.
\begin{theorem}
~
\begin{itemize}
\item For every $H\in(\frac{1}{2},1)$ and for every $k=1,\ldots,q-2$ we have
\begin{equation}
\lim_{_{N\rightarrow\infty}}N^{(2-2H^{%
\prime})(q-k)}T_{2q-2k}=z_{k,H}Z^{(2q-2k,(2q-2k)(H^{\prime}-1)+1)},\hskip%
0.5cm\mbox{ in }L^{2}(\Omega)  \label{convt2k}
\end{equation}
where $Z^{(2q-2k,(2q-2k)(H^{\prime}-1)+1)}$ denotes a Hermite random variable with self-similarity parameter
$(2q-2k)(H^{\prime}-1)+1$ and \\ $ z_{k,H}=d(H,q)^{2}a(H^{\prime})^{k}\left( (H^{\prime}-1)k+1\right) ^{-1}\left(
2(H^{\prime}-1)+1\right) ^{-1}.$
\item Moreover, if $H\in(\frac{3}{4},1)$ then
\begin{equation}
\lim_{_{N\rightarrow\infty}}N^{2-2H}x_{2,H}^{-1/2}T_{2q}=Z^{(2q,2H-1)},\hskip%
0.5cm\mbox{ in }L^{2}(\Omega).  \label{convfin}
\end{equation}
\end{itemize}
\end{theorem}
\textbf{Proof: } Recall that we have $%
T_{2k}=N^{2H-1}I_{2q-2k}\left(
\sum_{i=0}^{N-1}f_{i,N}\otimes_{k}f_{i,N}\right)$, for $k=1$ to $2q-2$, with $f_{i,N}%
\otimes_{k}f_{i,N}$ given by (\ref{fik}). The first step of the proof is to
observe that the limit of $N^{(2H^{\prime}-2)(q-k)}T_{2q-2k}$ is given by
\begin{equation*}
N^{(2H^{\prime}-2)(q-k)}N^{2H-1}I_{2q-2k}(f_{2q-2k}^{N})
\end{equation*}
with $ f_{2q-2k}^{N}(y_{1},\ldots,y_{q-k},z_{1},\ldots,z_{q-k})=d(H,q)^{2}a(H^{
\prime})^{k}\sum_{i=0}^{N-1}1_{[0,\frac{i}{N}]}(y_{i})1_{[0,\frac{i}{N} ]}(z_{i})\\
\int_{I_{i}}\int_{I_{i}}\partial_{1}K(u,y_{1})\ldots\partial_{1}K(u,y_{q-k})
\partial_{1}K(v,z_{1})\ldots\partial_{1}K(v,z_{q-k})|u-v|^{(2H^{
\prime}-2)k}dvdu.
$\newline
The argument leading to the above fact is the same as in Theorem \ref%
{ThmConv}: a the remainder term $r_{2q-2k}$, which is defined by $T_{2q-2k}$
minus $I_{2}(f_{2q-2k}^{N})$ converges to zero similarly to the term $r_{2}$
in the proof of Theorem \ref{ThmConv} because of the appearance of the
indicator functions $1_{I_{i}}(y_{i})$ or $1_{I_{i}}(z_{i})$ in each of the
terms that form this remainder.
The second step of the proof is to replace $\partial _{1}K(u,y_{i})$ by $%
\partial _{1}K(\frac{i}{N},y_{i})$ and $\partial _{1}K(v,z_{i})$ by $%
\partial _{1}K(\frac{i}{N},z_{i})$ on the interval $I_{i}$ inside the
integrals $du$ and $dv$. This can be argued, as in the proof of Theorem \ref%
{ThmConv}, by a dominated convergence theorem. Therefore the term $%
N^{(2-2H^{\prime })(q-k)}T_{2q-2k}$ will have the same limit as
\begin{eqnarray*}
&&N^{(2-2H^{\prime })(q-k)}N^{2H-1}d(H,q)^{2}a(H^{\prime
})^{k}\sum_{i=0}^{N-1}1_{[0,\frac{i}{N}]}(y_{i})1_{[0,\frac{i}{N}]}(z_{i}) \\
&&\prod_{j=1}^{q-k}\left[\partial _{1}K(\frac{i}{N},y_{j})\partial _{1}K(\frac{i}{N},z_{j})\right]\iint_{I_{i}^{2}}|u-v|^{(2-2H^{\prime })k}dvdu \\
&=&\frac{d(H,q)^{2}a(H^{\prime })^{k}}{\left( (H^{\prime }-1)k+1\right)
\left( 2(H^{\prime }-1)+1\right) }N^{(2-2H^{\prime
})(q-k)}N^{2H-1}N^{(2-2H^{\prime })k+2} \\
&&\sum_{i=0}^{N-1}1_{[0,\frac{i}{N}]^{2}}(y_{i},z_{i})\prod_{j=1}^{q-k}\left[\partial _{1}K(\frac{i}{N},y_{j})\partial _{1}K(\frac{i}{N},z_{j})\right] \\
&=&\frac{z_{k,H}}{N}\sum_{i=0}^{N-1}1_{[0,\frac{i}{N}]^{2}}(y_{i},z_{i})%
\prod_{j=1}^{q-k}\left[\partial _{1}K(\frac{i}{N},y_{j})\partial _{1}K(\frac{i}{N},z_{j})\right].
\end{eqnarray*}%
Now, the last sum is a Riemann sum that converges pointwise and in $ L^{2}([0,1]^{2q-2k})$ to $z_{k,H}L_{1}$ where
\begin{equation*}
L_{1}(y_{1},\ldots ,y_{q-k},z_{1},\ldots ,z_{q-k})=\int_{y_{1}\vee
z_{q-k}}^{1}\prod_{j=1}^{q-k}\left[\partial _{1}K(u,y_{j})\partial _{1}K(u,z_{j})\right]du
\end{equation*}
which is the kernel of the Hermite random variable of order $2q-2k$ with
self-similarity parameter $(2q-2k)(H^{\prime }-1)+1$. The case $H\in (\frac{3}{4}%
,1)$ and $k=0$ follows analogously. \hfill \vrule width.25cm height.25cm
depth0cm\smallskip

\section{Consistent estimation of $H$ and its asymptotics\label{EstimSec}}
Theorem \ref{ThmConv} can be immediately applied to the statistical
estimation of $H$. We note that with $V_{N}$ in (\ref{VN})\ and $S_{N}$
defined by%
\begin{equation}
S_{N}=\frac{1}{N}\sum_{i=0}^{N-1}\left( Z_{\frac{i+1}{N}}^{(q)}-Z_{\frac {i}{%
N}}^{(q)}\right) ^{2},  \label{SN}
\end{equation}
we have%
\begin{equation}
1+V_{N}=N^{2H}S_{N}  \label{VNSN}
\end{equation}
and $ \mathbf{E}\left[ S_{N}\right] =N^{-2H}$ so that $ H=-\frac{\log\mathbf{E}\left[ S_{N}\right] }{2\log N}.$  To form an
estimator of $H$, since Theorem \ref{ThmConv} implies that $S_{N}$ evidently concentrates around its mean, we will use $S_{N}$
in the role of
an empirical mean, instead of its true mean; in other words we let%
\begin{equation*}
\hat{H}=\hat{H}_{N}=-\frac{\log S_{N}}{2\log N}.
\end{equation*}
We immediately get from (\ref{VNSN}) that
\begin{eqnarray}
\log\left( 1+V_{N}\right) & =2H\log N+\log S_{N}   =2\left( H-\hat{H}\right) \log N.  \label{log1+VN}
\end{eqnarray}
The first observation is that $\hat{H}_{N}$ is strongly consistent for the
Hurst parameter.
\begin{proposition}
We have that $\hat{H}_{N}$ converges to $H$ almost surely as $N\rightarrow
\infty$.
\end{proposition}
\textbf{Proof: }Let us prove that $V_{N}$ converges to zero almost surely as $N\rightarrow\infty$. We know that $V_{N}$
converges to 0 in $L^{2}(\Omega)$ as $N\rightarrow\infty$ and an estimation for its variance is given by the formula
(\ref{meanT2}). On the other hand this sequence is stationary because the increments of the Hermite process are stationary. We
can therefore apply a standard argument to obtain the almost sure convergence for discrete stationary sequence under condition
(\ref{meanT2}). This argument follows from Theorem 6.2, page 492 in \cite{doob} and it can be used exactly as in the proof of
Proposition 1 in \cite{coeur}. A direct elementary proof using the Borel-Cantelli lemma is almost as easy. The almost sure
convergence of $\hat{H}$ to $H$ is then obtained immediately via (\ref{log1+VN}). \hfill\vrule width.25cm height.25cm
depth0cm\smallskip
\vskip0.5cm
Owing to the fact that $V_{N}$ is of small magnitude (it converges to zero
almost surely by the last proposition's proof), $\log\left( 1+V_{N}\right) $
can be confused with $V_{N}$. It then stands to reason that $\left( H-\hat
{%
H}\right) N^{2-2H^{\prime}}\log N$ is asymptotically Rosenblatt-distributed
since the same holds for $V_{N}$ by Theorem \ref{ThmConv}. Just as in that
theorem, more is true: as we now show, the convergence to a Rosenblatt
random variable occurs in $L^{2}\left( \Omega\right) $.
\begin{proposition}
\label{propfirstL2}There is a standard Rosenblatt random variable $R$ with
self-similarity parameter $2H^{\prime}-1$ such that%
\begin{equation*}
\lim_{N\rightarrow\infty}\mathbf{E}\left[ \left\vert
2N^{2-2H^{\prime}}\left( H-\hat{H}\right) \log
N-c_{2}c_{1,H,q}^{1/2}R\right\vert ^{2}\right] =0
\end{equation*}
\end{proposition}
\textbf{Proof: }To simplify the notation, we denote $c_{2}c_{1,H,q}^{1/2}$
by $c$ in this proof. Theorem \ref{ThmConv} signifies that a standard
Rosenblatt r.v. $R$ with parameter $2H^{\prime }-1$ exists such that%
\begin{equation*}
\lim_{N\rightarrow \infty }\mathbf{E}\left[ \left\vert N^{2-2H^{\prime
}}V_{N}-cR\right\vert ^{2}\right] =0.
\end{equation*}%
From (\ref{log1+VN}) we immediately get%
\begin{eqnarray*}
&&\mathbf{E}\left[ \left\vert 2N^{2-2H^{\prime }}\left( H-\hat{H}\right) \log N-cR\right\vert ^{2}\right] =\mathbf{E}\left[
\left\vert N^{2-2H^{\prime }}\log \left( 1+V_{N}\right)
-cR\right\vert ^{2}\right]  \\
&\leq &2\mathbf{E}\left[ \left\vert N^{2-2H^{\prime }}V_{N}-cR\right\vert
^{2}\right] +2N^{4-4H^{\prime }}\mathbf{E}\left[ \left\vert V_{N}-\log
\left( 1+V_{N}\right) \right\vert ^{2}\right] .
\end{eqnarray*}%
Therefore, we only need to show that $\mathbf{E}\left[ \left\vert V_{N}-\log
\left( 1+V_{N}\right) \right\vert ^{2}\right] =o\left( N^{4H^{\prime
}-4}\right) $.
Using the inequality $x-\log (1+x)=\vert x -\log (1-x) \vert \leq x^{2}$ for $x\geq \frac {1}{2}$ one gets
\begin{eqnarray*}
&&\mathbf{E}\left| V_{N} -\log (1+ V_{N}) \right| ^{2} \\
&&\leq 2 \mathbf{E} \left| V_{N} -\log (1+ V_{N}) \right| ^{2} 1_{V_{N}\geq -\frac{1}{2}}+ 2\mathbf{E} \left| V_{N} -\log (1+ V_{N}) \right| ^{2} 1_{V_{N}< -\frac{1}{2}}\\
&&\leq 2 \mathbf{E}V_{N}^{4} +4 P\left( V_{N} <-\frac{1}{2}\right) +4 \mathbf{E} \left| \log (1+V_{N})\right| ^{2} 1_{V_{N}< -\frac{1}{2}}.
\end{eqnarray*}
The first term  is bounded above and this is immediately dealt with using the
following lemma.
\begin{lemma}
\label{lemVN4}For every $n\geq2$, there is a constant $c_{n}$ such that $%
\mathbf{E}\left[ \left\vert V_{N}\right\vert ^{2n}\right] \leq c_{n}\ N^{(4H^{\prime}-4)n}.$
\end{lemma}
\textbf{Proof: } In the proof of Proposition \ref{PropLim} we calculated the
$L^{2}$ norm of the terms appearing in the decomposition of $V_{N}$, where $%
V_{N}$ is a sum of multiple integrals. Therefore, from Proposition 5.1 of
\cite{M} we immediately get an estimate for any event moment of each term.
Indeed, for $k=q-1$
\begin{equation*}
\mathbf{E}\left[ T_{2}^{2n}\right] \leq 1\cdot 3\cdot 5\cdots (4n-1)\;\big(%
c_{1,H}^{2}\;N^{4H^{\prime }-4}\big)^{n}=c_{n}N^{(4H^{\prime }-4)n}
\end{equation*}%
and for $1\leq k<q-1$
\begin{eqnarray*}
\mathbf{E}\left[ T_{2q-2k}^{2n}\right]  &\leq &1\cdot 3\cdot 5\cdots
(2(2q-2k)n)\;\big(c_{1,q,H}^{2}\;N^{(2q-2k)(2H^{\prime }-2)}\big)^{n} \\
&=&c_{q,n}N^{(2q-2k)(2H^{\prime }-2)n}
\end{eqnarray*}%
The dominant term is again the $T_{2}$ term and the result follows immediately.\hfill \vrule width.25cm height.25cm
depth0cm\smallskip \vskip0.5cm Then applying the above lemma for $n=2$ and using a similar proof to  prove that  $\mathbf{E}
\left| \log (1+V_{N})\right| ^{2} 1_{V_{N}< -\frac{1}{2}}=o(N^{4H'-4})$ and $P(V_{N}
<-\frac{1}{2})=O(N^{8H'-8})=o(N^{4H'-4})$. This leads to Proposition \ref{propfirstL2} \qed

A difficulty arises when applying the above proposition for model validation
when checking the asymptotic distribution of the estimator $\hat{H}$: the
normalization constant $N^{2-2H^{\prime}}\log N$ depends on $H$. While it is
not always obvious that one may replace this instance of $H^{\prime}$ by its
estimator, in our situation, because of the $L^{2}\left( \Omega\right) $
convergences, this is legitimate, as the following theorem shows.
\begin{theorem}
\label{ThmEstim}Let $\hat{H}^{\prime}=1+(\hat{H}-1)/q$. There is a standard
Rosenblatt random variable $R$ with self-similarity parameter $2H^{\prime}-1
$ such that%
\begin{equation*}
\lim_{N\rightarrow\infty}\mathbf{E}\left[ \left\vert 2N^{2-2\hat{H}%
^{\prime}}\left( H-\hat{H}\right) \log N-c_{2}c_{1,H,q}^{1/2}R\right\vert %
\right] =0.
\end{equation*}
\end{theorem}
\textbf{Proof: }By the previous proposition, it is sufficient to prove that%
\begin{equation*}
\lim_{N\rightarrow \infty }\mathbf{E}\left[ \left\vert \left( N^{2-2\hat{H}%
^{\prime }}-N^{2-2H^{\prime }}\right) \left( H-\hat{H}\right) \log
N\right\vert \right] =0.
\end{equation*}%
We decompose the probability space depending on whether $\hat{H}$ is far or
not from its mean. For a fixed value $\varepsilon >0$ which will be chosen
later, it is most convenient to define the event%
\begin{equation*}
D=\left\{ \hat{H}>q\varepsilon /2+2H-1\right\} .
\end{equation*}%
We have%
\begin{eqnarray*}
&&\mathbf{E}\left[ \left\vert \left( N^{2-2\hat{H}^{\prime
}}-N^{2-2H^{\prime }}\right) \left( H-\hat{H}\right) \log N\right\vert %
\right]  \\
&=&\mathbf{E}\left[ \mathbf{1}_{D}\left\vert \left( N^{2-2\hat{H}^{\prime
}}-N^{2-2H^{\prime }}\right) \left( H-\hat{H}\right) \log N\right\vert %
\right]  \\
&&+\mathbf{E}\left[ \mathbf{1}_{D^{c}}\left\vert \left( N^{2-2\hat{H}%
^{\prime }}-N^{2-2H^{\prime }}\right) \left( H-\hat{H}\right) \log
N\right\vert \right] =:A+B.
\end{eqnarray*}
We study $A$ first. Introducing the shorthand notation $x=\max \left(
2-2H^{\prime },2-2\hat{H}^{\prime }\right) $ and $y=\min \left( 2-2H^{\prime
},2-2\hat{H}^{\prime }\right) $, we may write
\begin{eqnarray*}
&& \left\vert N^{2-2\hat{H}^{\prime }}-N^{2-2H^{\prime }}\right\vert
=e^{x\log N}-e^{y\log N} \\
&=&e^{y\log N}\left( e^{\left( x-y\right) \log N}-1\right)
\leq  N^{y}\left( \log N\right) \left( x-y\right) N^{x-y} \\
&=&2\left( \log N\right) N^{x}\left\vert H^{\prime }-\hat{H}^{\prime
}\right\vert  = 2q^{-1}\left( \log N\right) N^{x}\left\vert H-\hat{H}\right\vert .
\end{eqnarray*}%
Thus%
\begin{eqnarray*}
A &\leq &2q^{-1}~\mathbf{E}\left[ \mathbf{1}_{D}N^{x}\left\vert H-\hat{H}%
\right\vert ^{2}\log ^{2}N\right]  \\
&=&2q^{-1}~\mathbf{E}\left[ N^{x-\left( 4-4H^{\prime }\right) }\mathbf{1}%
_{D}N^{4-4H^{\prime }}\left\vert H-\hat{H}\right\vert ^{2}\log ^{2}N\right] .
\end{eqnarray*}%
Now choose $\varepsilon \in (0,2-2H^{\prime })$. In this case, if $\omega
\in D$, and $x=2-2H^{\prime }$, we get $x-\left( 4-4H^{\prime }\right)
=-x<-\varepsilon $. On the other hand, for $\omega \in D$ and $x=2-2\hat{H}%
^{\prime }$, we get $x-\left( 4-4H^{\prime }\right) =2-2\hat{H}^{\prime
}-(4-4H^{\prime })<-\varepsilon $ as well. In conclusion, on $D$,%
\begin{equation*}
x-\left( 4-4H^{\prime }\right) <-\varepsilon ,
\end{equation*}%
which implies immediately%
\begin{equation*}
A\leq N^{-\varepsilon }2q^{-1}\mathbf{E}\left[ N^{4-4H^{\prime }}\left\vert
H-\hat{H}\right\vert ^{2}\log ^{2}N\right] ;
\end{equation*}%
this prove that $A$ tends to $0$ as $N\rightarrow \infty $, since the last
expectation above is bounded (converges to a constant) by Proposition \ref%
{propfirstL2}.
Now we may study $B$. We are now operating with $\omega \in D^{c}$. In other
words,
\begin{equation*}
H-\hat{H}>1-H-q\varepsilon /2.
\end{equation*}%
Still using $\varepsilon <2-2H$', this implies that $H>\hat{H}$.
Consequently, it is not inefficient to bound $\left\vert N^{2-2H^{\prime
}}-N^{2-2\hat{H}^{\prime }}\right\vert $ above by $N^{2-2\hat{H}^{\prime }}$%
. In the same fashion, we bound $\left\vert H-\hat{H}\right\vert $ above by $%
H$. Hence we have, using H\"{o}lder's inequality with the powers $2q$ and $%
p^{-1}+\left( 2q\right) ^{-1}=1$.%
\begin{eqnarray}
B &\leq &H\log N~\mathbf{E}\left[ \mathbf{1}_{D^{c}}N^{2-2\hat{H}^{\prime }}%
\right]   \label{BUL} \\
&\leq &H\log N~\mathbf{P}^{1/p}\left[ D^{c}\right] ~\mathbf{E}^{1/\left(
2q\right) }\left[ N^{\left( 2-2\hat{H}^{\prime }\right) 2q}\right] .  \notag
\end{eqnarray}%
From Proposition \ref{propfirstL2}, by Chebyshev's inequality, we have%
\begin{equation}
\mathbf{P}^{1/p}\left[ D^{c}\right] \leq \frac{\mathbf{E}^{1/p}\left[
\left\vert H-\hat{H}\right\vert ^{2}\right] }{\left( 1-H-q\varepsilon
/2\right) ^{2/p}}\leq c_{q,H}N^{-\left( 4-4H^{\prime }\right) /p}
\label{forBUL1}
\end{equation}%
for some constant $c_{q,H}$ depending only on $q$ and $H$. Dealing with the
other term in the upper bound for $B$ is a little less obvious. We must
return to the definition of $\hat{H}$. By (\ref{log1+VN}) we have%
\begin{eqnarray*}
1+V_{N} &=&N^{2\left( H-\hat{H}\right) }=N^{2q\left( H^{\prime }-\hat{H}%
^{\prime }\right) } \\
&=&N^{2q\left( 2-2\hat{H}^{\prime }\right) }N^{-2q\left( 2-2H^{\prime
}\right) }.
\end{eqnarray*}%
Therefore,
\begin{eqnarray}
\mathbf{E}^{1/\left( 2q\right) }\left[ N^{\left( 2-2\hat{H}^{\prime }\right)
2q}\right]  &\leq &N^{2-2H^{\prime }}\mathbf{E}^{1/\left( 2q\right) }\left[
1+V_{N}\right]   \notag \\
&\leq &2N^{2-2H^{\prime }}.  \label{forBUL2}
\end{eqnarray}%
Plugging (\ref{forBUL1}) and (\ref{forBUL2}) back into (\ref{BUL}), we get%
\begin{equation*}
B\leq 2Hc_{q,H}\left( \log N\right) N^{-\left( 4-4H^{\prime }\right) \left(
2/p-1\right) }.
\end{equation*}%
Given that $2q\geq 4$ and $p$ is conjugate to $2q$, we have $p\leq 4/3$ so
that $2/p-1\geq 1/2>0$, which proves that $B$ goes to $0$ as $N\rightarrow
\infty $. This finishes the proof of the theorem.\hfill \vrule width.25cm
height.25cm depth0cm\bigskip

Finally we state the extension of our results to $L^{p}\left( \Omega\right) $%
-convergence.
\begin{theorem}
The convergence in Theorem \ref{ThmEstim} holds in any $L^{p}\left(
\Omega\right) $. In fact, the $L^{2}\left( \Omega\right) $-convergences in
all other results in this paper can be replaced by $L^{p}\left(
\Omega\right) $-convergences.
\end{theorem}
\textbf{Sketch of proof}. We only give the outline of the proof. Lemma \ref%
{lemVN4} works because, in analogy to the Gaussian case, for random
variables in a fixed Wiener chaos of order $p$, existence of second moments
implies existence of all moments, and the relations between the various
moments are given using a set of constants which depend only on $p$. This
Lemma can be used immediately to prove the extension of Proposition \ref%
{PropLim} that
\begin{equation*}
\mathbf{E}\left[ V_{N}^{2p}\right] \simeq c_{p}N^{p(4H^{\prime} -4)}.
\end{equation*}%
Proving a new version of Theorem \ref{ThmConv} with $L^{p}\left( \Omega
\right) $-convergence can base itself on the above result, and requires a
careful reevaluation of the various terms involved. Proposition \ref%
{propfirstL2} can then be extended to $L^{p}\left( \Omega \right) $%
-convergence thanks to the new $L^{p}\left( \Omega \right) $ versions of
Theorem \ref{ThmConv} and Proposition \ref{PropLim}, and Theorem \ref%
{ThmEstim} follows easily from this new version of Proposition \ref{PropLim}%
.\qed


\begin{thebibliography}{99}
\bibitem{Beran} {J. Beran (1994): }\emph{Statistics for Long-Memory
Processes. } {Chapman and Hall. }
\bibitem{BrMa} {P. Breuer and P. Major (1983): }Central limit theorems for
nonlinear functionals of Gaussian fields. \emph{J. Multivariate Analysis}{,
\textbf{13 }(3), 425-441. }
\bibitem{BuSw}
{K. Burdzy and J. Swanson (2008): }Variation for the solution to the stochastic heat equation II. Preprint.

\bibitem{coeur} {J.F. Coeurjolly (2001): }Estimating the parameters of a
fractional Brownian motion by discrete variations of its sample paths. \emph{%
Statistical Inference for Stochastic Processes}{, \textbf{4}, 199-227. }
\bibitem{DM} {R.L. Dobrushin and P. Major (1979): }Non-central limit
theorems for non-linear functionals of Gaussian fields. \emph{Z.
Wahrscheinlichkeitstheorie verw. Gebiete}{, \textbf{50}, 27-52. }
\bibitem{doob} {J.L. Doob (1953): }\emph{Stochastic Processes. }{Wiley
Classics Library. }
\bibitem{GS} L. Giraitis and D. Surgailis (1985): \textit{CLT and other
limit theorems for functionals of Gaussian processes.} Z. Wahrsch. verw.
Gebiete \textbf{70}, 191-212.
\bibitem{GuLe} {X. Guyon and J. Le\'{o}n (1989): }Convergence en loi des $H$%
-variations d'un processus gaussien stationnaire sur $\mathbf{R}$. \emph{%
Annales IHP}{, \textbf{25}, 265-282. }
\bibitem{LaIs} {J. Istas and G. Lang  (1997): }Quadratic variations and
estimators of the H\"{o}lder index of a Gaussian process. \emph{Annales IHP}{%
, \textbf{33}, 407-436. }
\bibitem{M} {P. Major (2005): }\emph{Tail behavior of multiple random
integrals and $U$-statistics. }{Probability Surveys.}

\bibitem {No}{I. Nourdin (2008): } Asymptotic behavior of certain weighted
quadratic variations and cubic varitions of fractional Brownian motion. \emph{The Annals of Probability, }{\textbf{6}, 2159-2175. }
\bibitem{NoPe}{I. Nourdin and G. Peccati (2008): }Weighted power variations of iterated Brownian motion. \emph{Electronic Journal of Probability, } {\textbf{13} (43), 1229-1256. }

\bibitem{NNT} {I. Nourdin, D. Nualart and C.A. Tudor (2007): } \emph{Central
and non-central limit theorems for weighted power variations of fractional
Brownian motion. }\emph{Annales IHP, } to appear.
\bibitem{N} {D. Nualart (2006): }\emph{Malliavin Calculus and Related
Topics. Second Edition. }{Springer. }

\bibitem{Rev} {A. R\'eveillac (2009): }Convergence of finite-dimensional laws of the weighted quadratic variations process for some
 fractional Brownian sheets.  \emph{Stoch. Anal. Appl., } {\textbf{27},  no. 1, 51--73.}
\bibitem {Sw}{J. Swanson (2007): }Variations of the solution to a stochastic
heat equation. \emph{Annals of Probability}{, \textbf{35}(6), 2122--2159.}

\bibitem{Ta1} {M. Taqqu (1975)}: Weak convergence to the fractional Brownian
motion and to the Rosenblatt process. \emph{Z. Wahrscheinlichkeitstheorie
verw. Gebiete}{, \textbf{31}, 287-302.}
\bibitem{Ta2} {M. Taqqu (1979}): Convergence of integrated processes of
arbitrary Hermite rank. \emph{Z. Wahrscheinlichkeitstheorie verw. Gebiete}{,
\textbf{50}, 53-83. }
\bibitem{T}
{C.A. Tudor (2008): }Analysis of the Rosenblatt process. \emph{ESAIM Probability and Statistics}{, \textbf{12}, 230-257.}
\bibitem{TV} {C.A. Tudor and F. Viens (2009): }\emph{Variations and
estimators for the selfsimilarity order through Malliavin calculus. }\emph{ The Annals of Probability, } {\textbf{6}, 2093-2134. }
\end{thebibliography}
\end{document}